\documentclass[final,5p,times,twocolumn]{elsarticle}
\usepackage{etex}
\usepackage{lineno,hyperref}
\usepackage{amsmath}
\usepackage{amsfonts}
\modulolinenumbers[5]
\journal{CMAME}

\setcitestyle{square}
\biboptions{numbers,sort&compress}
\modulolinenumbers[5]
\usepackage{multirow}
\usepackage{array}
\usepackage{subfig}
\usepackage{gensymb}
\usepackage{cancel}
\usepackage{xcolor}
\usepackage{empheq}
\usepackage[many]{tcolorbox}
\usepackage{titlesec}
\titlelabel{\thetitle.\ }
\usepackage[symbol]{footmisc}
\usepackage[linesnumbered,ruled,vlined]{algorithm2e}
\usepackage{morefloats}
\usepackage{mathrsfs}
\usepackage{relsize}
\usepackage{todonotes}

\usepackage{breqn}
\usepackage{textcomp}
\DeclareMathAlphabet\mathbfcal{OMS}{cmsy}{b}{n}
\usepackage[many]{tcolorbox}
\newtcolorbox{highlighted}{colback=yellow,coltext=black,breakable}
\makeatletter
\newcommand{\mathleft}{\@fleqntrue\@mathmargin\parindent}
\newcommand{\mathcenter}{\@fleqnfalse}
\makeatother
\usepackage{acronym}
\usepackage[autostyle]{csquotes}
\usepackage[figurename=Fig.]{caption}

\let\orgautoref\autoref

\renewcommand{\autoref}[1]
{%
\def\equationautorefname{Eq.}%
\def\figureautorefname{Fig.}%
\def\subfigureautorefname{Fig.}%
\orgautoref{#1}%
}

\acrodef{FEA}{Finite Element Analysis}
\acrodef{ECP}{Element Connectivity Parameterization}
\acrodef{ESO}{Evolutionary Structural Optimization}
\acrodef{DOF}{Degree-Of-Freedom}
\acrodef{DOFs}{Degrees-Of-Freedom}
\acrodef{HPM}{Heaviside Projection Method}
\acrodef{BESO}{Bidirectional Evolutionary Structural Optimization}
\acrodef{NAND}{Nested ANalysis and Design}
\acrodef{SIMP}{Solid Isotropic Material with Penalization}
\acrodef{RAMP}{Rational Approximation of Material Properties}
\acrodef{MMA}{Method of Moving Asymptotes}
\acrodef{GMRES}{Generalized Minimal RESidual}
\acrodef{TO}{Topology Optimization}
\acrodef{OC}{Optimality Criteria}
\acrodef{PDE}{Partial Differential Equation}
\acrodef{GEP}{Generalized Eigenvalue Problem}
\acrodef{BLF}{Buckling Load Factor}
\acrodef{IG}{Initial Guess}

\graphicspath{ {./figs_new/} }

\makeatletter
\def\ps@pprintTitle{%
	\let\@oddhead\@empty
	\let\@evenhead\@empty
	\def\@oddfoot{\footnotesize\itshape
		{} \hfill\today}%
	\let\@evenfoot\@oddfoot
}
\makeatother

\begin{document}

\begin{frontmatter}

\title{Revisiting element removal for density-based structural topology optimization with reintroduction by Heaviside projection}

\author[JHU]{Reza Behrou\corref{cor1}}
\ead{rbehrou@jhu.edu}
\author[JHU]{Reza Lotfi}
\author[JHU,MIT]{Josephine Voigt Carstensen}
\author[JHU]{Federico Ferrari}
\author[JHU]{James K. Guest}

\cortext[cor1]{Corresponding Author: Reza Behrou}

\address[JHU]{Department of Civil and Systems Engineering, Johns Hopkins University, Baltimore, MD 21218, USA}
\address[MIT]{Department of Civil Engineering, Massachusetts Institute of Technology, Cambridge, MA 02139, USA}

%
\begin{abstract}
We present a strategy grounded in the element removal idea of Bruns and Tortorelli \cite{BT:03} and aimed at reducing computational cost and circumventing potential numerical instabilities of density-based topology optimization. The design variables and the  relative densities are both represented on a fixed, uniform finite element grid, and linked through filtering and Heaviside projection. The regions in the analysis domain where the relative density is below a specified threshold are removed from the forward analysis and replaced by fictitious nodal boundary conditions. This brings a progressive cut of the computational cost as the optimization proceeds and helps to mitigate numerical instabilities associated with low-density regions. Removed regions can be readily reintroduced since all the design variables remain active and are modeled in the formal sensitivity analysis. A key feature of the proposed approach is that the Heaviside functions promote material reintroduction along the structural boundaries by amplifying the magnitude of the sensitivities inside the filter reach. Several 2D and 3D structural topology optimization examples are presented, including linear and nonlinear compliance minimization, the design of a force inverter, and frequency and buckling load maximization. The approach is shown to be effective at producing optimized designs equivalent or nearly equivalent to those obtained without the element removal, while providing remarkable computational savings.
\end{abstract}


\end{frontmatter}
%

\section{Introduction}
 \label{sec:Introduction}

This paper presents an adaptive strategy for reducing the computational cost of density-based topology optimization. The strategy systematically removes elements whose relative density is below a user-defined threshold from the \ac{FEA} (i.e., reducing the dimensionality of the analysis space) without precluding their reappearance (i.e., preserving the freedom of the original design space).

In the classical approach of density-based topology optimization, the entire design domain is modeled in both the analysis and re-design steps, including the regions where the relative density reaches its lower bound (e.g., voids). These \enquote{void} regions are clearly unnecessary, neither for the physical simulation, as they bring negligible contributions to the system's response, nor for manufacturing purposes, as they are replaced by holes during the design post-processing. Their presence in the computational domain basically serves only to allow the reintroduction of material in the optimization procedures. However, low-density regions are quite detrimental to the analysis, as they (1) increase the size of the system to be solved, with resulting computational burden, and (2) are prone to numerical instabilities, especially when nonlinear physics are considered.

Reducing the computational cost of topology optimization procedures is a key requirement for both industrial and academic purposes and a pivotal point in the strive to make topology optimization a comprehensive and fast design tool \cite{DG:14}. In the nested approach to topology optimization, the computational overhead is mainly driven by the \ac{FEA} steps (i.e., assembly and solving). Researchers have tackled this issue mainly by employing efficient reanalysis techniques \cite{Kirsch:03,KB:04,ABS:09}, efficient preconditioned solvers \cite{WSP:07,ASS:10,AS:11,FLS:18,FS:20a}, remeshing strategies \cite{MR:95, Stainko:06, SPW:08, LHZ+:18}, or reduced-order models \cite{Yoon:10, XLB+:20}.

Circumventing numerical instabilities associated with low-density elements has drawn many researchers’ attention, especially for problems governed by nonlinear mechanics and eigenvalue analysis. In large deformation settings, low-density elements may be prone to distortion and/or prevent convergence of the nonlinear analysis. Researchers have tackled this in different ways, such as removing the \ac{DOFs} associated with low-stiffness elements from the equilibrium convergence criteria \cite{BPS:00, PBS:01}, stiffening elements that undergo large deformations using a hyperelastic constitutive model \cite{BT:01, KS:13}, using stiff continuum elements connected by rigid links that vanish at void locations \cite{YK:07}, and modifying the interpolation of strain energy to be linear in low-density elements \cite{WLS+:14}. For eigenvalue problems, artificial (buckling and vibration) modes appearing in zones with low relative density are classically cured by modifying the interpolation function \cite{Pedersen:00,DG:14} in such regions.

Here we revisit and enhance the approach proposed by \cite{BT:03}: we remove elements with a relative density below the user-defined threshold value from the \ac{FEA}, and artificial nodal boundary conditions suppress the \ac{DOFs} surrounded by these removed elements. Therefore, the size of the finite element system to be solved is reduced during the optimization procedure, with a reduction ratio that is roughly proportional to the allowed volume fraction, considerably cutting the computational cost. Additionally, removing elements with low relative density from the \ac{FEA} naturally circumvents instabilities associated with them.

A key feature of the proposed approach is that material may be readily reintroduced into the domain of removed elements, freeing \ac{DOFs} that may have been artificially restrained. This is an important property, particularly for more complex topology optimization problems where structural features may translate across or re-appear in the design domain during the design evolution. The ability to reintroduce material is due to the separation of the design and analysis spaces, and governed purely by the formal sensitivity propagation due to the filter and projection in the \ac{HPM}. Although the original work \cite{BT:03} was likewise able to reintroduce material, this effect was limited due to using a linear density filter alone. We will show how the use of the \ac{HPM} \cite{GPB:04} greatly benefits the material reintroduction, making it faster and suitable for more complex problems.  Additionally, \cite{BT:03} reported that load paths could become disconnected using the element removal strategy, and therefore recommended delaying element removal until structural topology had formed. The \ac{HPM}-based approach herein uses standard tightened optimization parameters as described in \cite{GAH:11} and no instances of disconnected structures was observed, allowing the approach to be implemented from the first optimization iteration.   

We note that other works have proposed reducing the number of \ac{DOFs} in the system by using adaptive meshing strategies, either to re-mesh the entire domain (void and solid) \cite{Stainko:06, SPW:08} or to create discretizations that are body-fitted (to varying degree) to the current, evolving design \cite{MR:95, LHZ+:18}. Maute and Ramm \cite{MR:95}, for example, implemented a formal re-meshing scheme where relative density isolines were computed on the design space, and a threshold line was selected to represent the structural boundary. As this boundary evolves, the structural domain is re-meshed. This adaptive scheme was extended to problems governed by elastoplastic material models and was capable of circumventing instabilities in low stiffness elements \cite{MSR:98}, just as the proposed approach. Researchers have also taken advantage of the separation of the design and analysis meshes in projection by meshing the two spaces differently \cite{MR:95, GS:10, NPS+:10}, including with adaptive features \cite{MR:95, GS:10}. We finally remark that the single-phase modeling has long been a selling point of discrete implementations of level set methods \cite{WWG:03, AJT:04, Challis:10}, with savings specifically quantified in \cite{CG:09} for topology optimization of fluids. It has also been a promoted feature of hard-kill \ac{ESO} approaches \cite{XS:97, QSX:98}, where material may not return to elements once removed, as well as soft-kill methods such as \ac{BESO} where element reintroduction is driven by heuristic rules based on elemental parameters, such as strain energy \cite{ZR:01}. The proposed approach differs from these existing works in that (1) it operates on a fixed finite element mesh where only nodal boundary conditions are updated, and (2) material reintroduction is driven entirely by the gradient-based optimizer informed through formal sensitivity analysis of the current design. It is clearly acknowledged that introducing adaptive boundary conditions creates a discontinuity in the state variables of the optimization problem, but the findings presented herein demonstrate that the approach is remarkably robust for a variety of design problems.

The remainder of the paper is organized as follows. In Section \ref{sec:OptimizationProblem}, we outline the general formulation and ingredients of density-based topology optimization. Section \ref{sec:SensitivityStudy} illustrates the concepts at the core of the element removal strategy by comparing the influence of various filtering and projection methods on the sensitivity and topology propagation. The procedure implementing the adaptive elemental removal strategy, and the corresponding application of adaptive nodal boundary conditions is described in Section \ref{sec:NodalBoundaryCondition}. In Section \ref{sec:NumericalExamples} we apply the adaptive element removal method to several structural topology optimization problems, governed by linear, nonlinear and eigenvalue state equations, showing its ability to achieve optimized design nearly equivalent to those given by the \enquote{standard} approach, and with remarkable computational savings. Concluding remarks are summarized in Section \ref{sec:Conclusions}.

\section{Density-based topology optimization framework}
 \label{sec:OptimizationProblem}
We consider a structured discretization $\Omega_{h}$ consisting of $m$ equi-sized elements $\Omega_{e}$. The topology of a structure is parametrized by the (nodal- or element-based) design variables $\phi_{i}\in [0,1]$, $i=1,\ldots,n_{d}$. Let $\boldsymbol{u} \in \mathbb{R}^{n}$ represent the state of the system, described by $n$ \ac{DOFs}.

We seek the configuration $\boldsymbol{\phi}^{\ast}$ that minimizes an objective $g_{0}$, while a set of constraints $g_{j}$, $j=1,\ldots,n_{g}$ are simultaneously satisfied. As various topology optimization problem types are solved herein we express the optimization problem in general form as follows:
\begin{equation}
 \label{eq:OptimizationProblem_1}
  \begin{aligned}
  \min_{\boldsymbol{\phi}} \:
  & \quad g_{0}(\boldsymbol{\phi},\boldsymbol{u}(\boldsymbol{\phi})) \\
  \text{s.t.} &
  \begin{cases}
  \boldsymbol{r} ( \boldsymbol{\phi}, \boldsymbol{u}(\boldsymbol{\phi})) = \boldsymbol{0} & \\
  g_{j} (\boldsymbol{u}(\boldsymbol{\phi}), \boldsymbol{\phi})) \leq 0 & j = 1,\ldots, n_{g} \\ 
  0 \leq \phi_{i} \leq 1 & i = 1,\ldots,n_{d}
  \end{cases},
 \end{aligned}
\end{equation}
where $\boldsymbol{r} (\boldsymbol{\phi}, \boldsymbol{u}(\boldsymbol{\phi}))$ is the residual of the linear or nonlinear PDEs governing the state of the system.

To avoid mesh dependency, checkerboard patterns and other issues arising when the design variables $\boldsymbol{\phi}$ are directly used for determining the physical response and sensitivities \cite{SP:98}, the relative density variable $\rho_{e} = \rho_{e}(\boldsymbol{\phi})$ is introduced on each element $\Omega_{e}$, and related to the design variables by \ac{HPM} \cite{GPB:04}.

First, the intermediate variables $\mu_{e}(\boldsymbol{\phi})$, $i=1, \ldots, m$ are obtained through the filtering operation \cite{BT:01}:
\begin{equation}
 \label{eq:DO_SensitivityAnalysis_14}
 \mu_{e}(\boldsymbol{\phi}) = 
 \frac{\sum_{i\in\mathcal{N}_{i}}
 w(\boldsymbol{x}_{i}, \boldsymbol{x}_{e})\phi_{i}}
 {\sum_{i\in\mathcal{N}_{i}}
 w(\boldsymbol{x}_{i}, \boldsymbol{x}_{e})},
\end{equation}
where $w(\boldsymbol{x}_{i}, \boldsymbol{x}_{e})$ is the weighting kernel depending on the distance between the centroid location of element $\boldsymbol{x}_{e}$ and the location of the design variable $\boldsymbol{x}_{i}$, and $\mathcal{N}_{i}$ is the neighborhood set that contains all design variables located within a distance $r_{\rm min}$ of the centroid $e$. In this work, we adopt the hat weighting function
\begin{equation}
 \label{eq:filterHatFunction}
  w(\boldsymbol{x}_{i}, \boldsymbol{x}_{e}) =
  \begin{cases}
   1 - \frac{\|\boldsymbol{x}_{i} - \boldsymbol{x}_{e} \|_{2}}{r_{\rm min}} & {\rm if} \quad \|\boldsymbol{x}_{i} - \boldsymbol{x}_{e} \|_{2} \leq r_{\rm min} \\
   0 & {\rm if} \quad \|\boldsymbol{x}_{i} - \boldsymbol{x}_{e} \|_{2} > r_{\rm min}
  \end{cases},
\end{equation}
resulting in a linear relationship between $\boldsymbol{\phi}$ and $\boldsymbol{\mu}$.

On each element, $\rho_{e}$ is then related to $\mu_{e}$ by the relaxed Heaviside function \cite{GPB:04}
\begin{equation}
 \label{eq:DO_SensitivityAnalysis_15}
 \rho_{e} (\boldsymbol{\phi}) = 1 -
 e^{-\beta \mu_{e}(\boldsymbol{\phi})} +
 \mu_{e}(\boldsymbol{\phi})e^{-\beta},
\end{equation}
where $\beta \geq 0$ dictates the sharpness of the projection. As $\beta \to \infty$, \autoref{eq:DO_SensitivityAnalysis_15} approaches the step operator, achieving a 0-1 density distribution and imposing the minimum length scale on the solid features \cite{GPB:04, CG:18}.

The material properties of the system are then related to relative densities, $\rho_{e}$. Namely, given the physical property \enquote{$p$}, we consider a relationship in the form $p(\rho_{e},\eta,p_{\rm min},p_{0})$, where $p_{0}$ is the property value corresponding to the solid, $p_{\rm min}$ that of the (artificial) void material, and $\eta$ is a term that penalizes intermediate densities $\rho_{e}\in (0,1)$.
In this work, we will refer to some widely used interpolation schemes, such as the \ac{SIMP} \cite{Bendsoe:89} or the \ac{RAMP} \cite{SS:01}, and the particular interpolation adopted to represent physical properties such as stiffness, mass and stresses will be specified for each example in \autoref{sec:NumericalExamples}.

Assuming that the relationship $\boldsymbol{\phi}\mapsto \boldsymbol{u}(\boldsymbol{\phi})$ is well-defined (i.e., that the tangent operator $K := \partial_{\boldsymbol{u}} \boldsymbol{r} (\boldsymbol{\phi}, \boldsymbol{u}(\boldsymbol{\phi}))$ is non-singular on the feasible set), we solve \autoref{eq:OptimizationProblem_1} by a nested approach, consisting of sequential analysis and re-design steps \cite{BS:04}.

At each given design point $\boldsymbol{\phi}$, the state variables are computed by solving $\boldsymbol{r}(\boldsymbol{\phi}, \boldsymbol{u}(\boldsymbol{\phi})) = \boldsymbol{0}$; the objective and constraint functions and their sensitivities are evaluated; and the design variables are updated by solving a locally approximate optimization problem through a gradient-based method.

Given the non-singularity of $K(\boldsymbol{\phi})$, we may express
\begin{equation}
 \label{eq:DO_SensitivityAnalysis_1}
 \frac{d \boldsymbol{u}}{d \rho_{e}}(\boldsymbol{\phi}) = 
 - [K(\boldsymbol{\phi})]^{-1}
 \frac{\partial \boldsymbol{r}}
 {\partial \rho_{e}}(\boldsymbol{\phi},\boldsymbol{u}(\boldsymbol{\phi})),
\end{equation}
and the derivatives of the objective and constraint functions with respect to $\rho_{e}$ read
\begin{equation}
 \label{eq:DO_SensitivityAnalysis_2}
  \frac{d g_{j}}{d \rho_{e}}
  (\boldsymbol{\phi},\boldsymbol{u}(\boldsymbol{\phi})) =
  \frac{\partial g_{j}}{\partial \rho_{e}}
  (\boldsymbol{\phi},\boldsymbol{u}(\boldsymbol{\phi})) + \boldsymbol{a}^{T}_{j}
  \frac{\partial \boldsymbol{r}}{\partial \rho_{e}}
  (\boldsymbol{\phi},\boldsymbol{u}(\boldsymbol{\phi})),
\end{equation}
where $\boldsymbol{a}_{j}(\boldsymbol{\phi}) = -K^{-T}(\partial_{\boldsymbol{u}} g_{j}(\boldsymbol{\phi},\boldsymbol{u}(\boldsymbol{\phi})))$ is the adjoint vector that eliminates the implicit relationship $\boldsymbol{\phi} \mapsto \boldsymbol{u}(\boldsymbol{\phi})$.

Finally, the formal sensitivities with respect to the independent variables $\phi_{i}$ are recovered from \eqref{eq:DO_SensitivityAnalysis_2} by the chain rule:
\begin{equation}
 \label{eq:chainRule_Sensitivity}
  \frac{d g_{j}}{d \phi_{i}} =
  \sum_{e \in \mathcal{N}_{i}}
  \frac{d g_{j}}{d \rho_{e}}
  \frac{\partial \rho_{e}}{\partial \mu_{e}}
  \frac{\partial \mu_{e}}{\partial \phi_{i}},
\end{equation}
where $\mathcal{N}_{i}$ is the neighborhood set containing all elements within the distance $r_{\rm min}$ of the design variable $i$ and $\partial\rho_{e}/\partial\mu_{e} = \beta e^{-\beta \mu_{e}(\boldsymbol{\phi})}+e^{-\beta}$. The term $\partial\mu_{e}/\partial\phi_{i}$ involves the application of \autoref{eq:DO_SensitivityAnalysis_14}, with $\partial g_{j}/\partial\rho_{e}$ instead of $\boldsymbol{\phi}$.

\begin{figure*}[tb]
 \centering
  \includegraphics[width=0.98\linewidth]{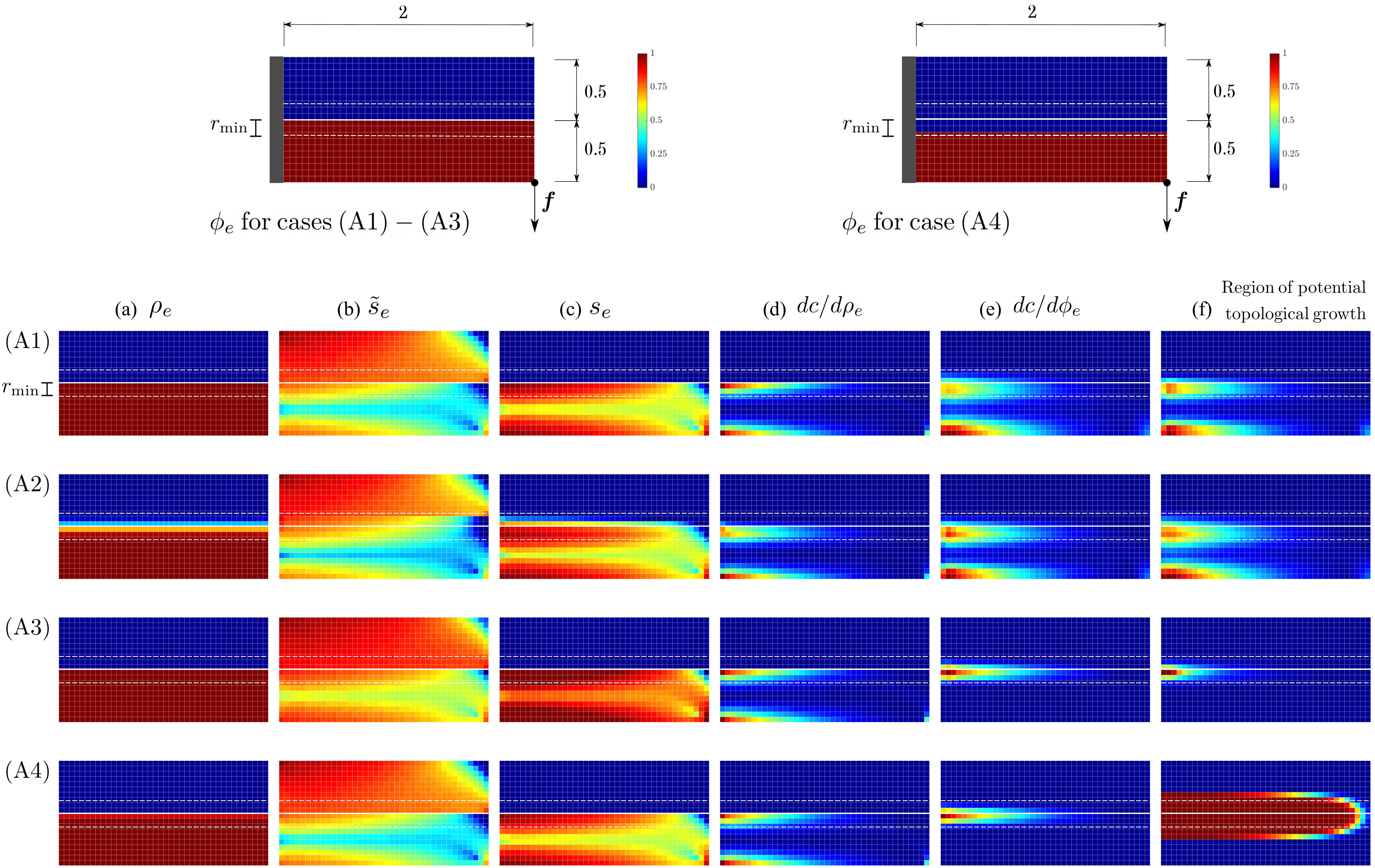}
\caption{Propagation of sensitivity and topology information when using different filtering and projection approaches: Sensitivity filter (A1) \cite{Sigmund:07}, Linear density filter (A2) \cite{BT:01}, threshold projection (A3) \cite{WLS:11}, and Heaviside projection method (A4) \cite{GPB:04}. The design variables are distributed according to the two sketches on top, such that $\rho_{e} = 1$ on the lower half of the domain and $\rho_{e} = 0$ in the upper half, for all the cases (some negligible differences along the separation line are due to the filter). The filter radius is set to $r_{\rm min} = 2.5$ elements for all cases, and $\beta = 40$ for both (A3) and (A4). The threshold value used for (A3) is 0.5. In plots of columns (a) and (f), the quantities range between 0 and 1, while plots in columns (b)-(e) show the $\log$-scale distribution of the quantities, normalized with respect to their maximum values.}
\label{fig:SensitivityStudy}
\end{figure*}

\section{Influence of \ac{HPM} on sensitivity and topology propagation}
 \label{sec:SensitivityStudy}

Before introducing the element removal strategy, we discuss the effect of the three main operations relating the design variables $\boldsymbol{\phi}$ to the physical response of the system. These are
\begin{itemize}
 \item Interpolation of the material properties having a penalization on the intermediate densities;
 \item Filtering (\autoref{eq:DO_SensitivityAnalysis_14}), that propagates information;
 \item Heaviside projection (\autoref{eq:DO_SensitivityAnalysis_15}), amplifying the information.
\end{itemize}

To this end, we consider the simple cantilever beam sketched in \autoref{fig:SensitivityStudy}, for the linearized compliance setting. The domain is discretized by $40 \times 20$ elements $\Omega_{e}$ and the filter radius is set to $r_{\rm min} = 2.5h_{e}$, where $h_{e}$ is the elemental length. We adopt the \ac{SIMP}-like interpolation for the Young's modulus $E(\rho_{e}) = \rho_{e}^{\eta} E_{0}$, where $E_{0} = 1$ and $\eta = 3$. Let's introduce the \textit{non-dimensional} elemental strain energy $\tilde{s}_{e} = (\boldsymbol{u}_{e})^{T}{K_{e}}_{0}\boldsymbol{u}_{e}$, where ${K_{e}}_{0}$ is the (non-dimensional) elemental stiffness matrix. The physically meaningful elemental strain energy is $s_{e} = E(\rho_{e})\tilde{s}_{e}$ and the response of interest for this example is the system compliance $c(\boldsymbol{\phi}) = \sum^{m}_{e=1} s_{e} \Omega_{e}$.

We compare four different approaches: (A1) Linear sensitivity filtering \cite{Sigmund:07}, (A2) Linear density filtering \cite{BT:01}, (A3) threshold projection \cite{WLS:11}, and (A4) \ac{HPM} \cite{GPB:04} (see rows in \autoref{fig:SensitivityStudy}). These introduce a different relationship between the design variables $\phi_{e}$ and the relative densities $\rho_{e}$ (both located at the element centroids here), and between the derivatives $dc/d\rho_{e}$ and $dc/d\phi_{e}$. In the following we formally denote these relationships by $\rho_{e} = \mathcal{T}(\phi_{i})$ and $dc/d\phi_{e} = \mathcal{T}^{-1}(dc/d\rho_{e})$, respectively.

For cases (A1)-(A2) we set $\phi_{e} = 1$ in the lower half of the beam and $\phi_{e} = 0$ in the upper half as shown in \autoref{fig:SensitivityStudy}. The corresponding volume fractions (i.e., $\rho_{e}$) are shown in column (a), where we note the density filter (A2) creates a blurry boundary, which also normally appears when using the sensitivity filter within topology optimization (though is not shown here as we have manually selected the $\phi_{e}$ distribution). The same $\phi$ distribution is used for thresholding projection (A3) with threshold equal to 0.5 and $\beta = 40$, leading to a nearly equivalent distribution of $\rho_{e}$. In order to fairly compare the reach of the sensitivities in all the approaches, for the \ac{HPM} case (A4) we set $\phi_{e} = 0$ for the top half of the beam and also on the first two element layers from the centerline in the lower half of the beam, which when filtered and projected (\autoref{eq:DO_SensitivityAnalysis_15}) with $\beta = 40$ creates the nearly identical topology $\rho_{e}$ shown in column (a)  (see \autoref{fig:SensitivityStudy}). 

The distributions of the non-dimensional $\tilde{s}$ and dimensional $s$ strain energies are shown in columns (b) and (c) of \autoref{fig:SensitivityStudy}, respectively. We note that each property is very similar across the four cases, with differences due to the small differences at the structural boundary seen in column (a). Another important observation is that, due to continuity of deformation,  $\tilde{s}$ may be large on elements with low relative density $\rho_{e}$ (top half of the beam). However, the \ac{SIMP} penalization drives the strain energy distribution $s_{e}$ to zero in these elements when $\eta > 1$. It then follows that the compliance derivatives $d c/d \rho_{e} = \eta\rho^{\eta-1}_{e}\tilde{s}_{e}$ are likewise driven to zero for $\eta > 2$, as shown in column (d) of \autoref{fig:SensitivityStudy}. Therefore, for a penalization exponent $\eta > 2$, elements $\Omega_{e}$ with low relative density (say, $\rho_{e} < 0.1$) carry a negligible contribution to both the physical response and derivatives $d c/d \rho_{e}$. This motivates the idea of removing them from the state analysis.

The corresponding derivatives of compliance with respect to the independent optimization variable, $d c/d \phi_{e} = \mathcal{T}^{-1}(d c/d \rho_{e})$, are shown in the column (e) of \autoref{fig:SensitivityStudy}. The key observation here is that if an element $\Omega_{e}$ is within the distance $r_{\rm min}$ from the structural boundary, its sensitivities $dc/d\phi_{e}$ can be non-zero even if its $d c/d\rho_{e} = 0$. This is purely due to the propagation effect under the filtering operation. This suggests that the material can be reintroduced on an element $\Omega_{e}$ that has been removed from the finite element analysis (thus, where $s_{e} = d c / d\rho_{e} = 0$), provided that one of its neighboring elements $\Omega_{j}$ within the filter reach has $d c/ d \rho_{j} \neq 0$.

Plots in column (f) display the \enquote{region of potential topological growth} for the four different approaches. These plots are obtained by applying the relationship $\rho_{e} = \mathcal{T}(d c/ d\phi_{e})$, thus treating the sensitivity coefficients in column (e) formally as design variables. This operation is not physically meaningful, as we might have $\rho_{e} \notin [0,1]$. Nevertheless, it clearly shows the topological influence of the two combined operations $dc/d \phi_{e} = \mathcal{T}^{-1}(d c/d \rho_{e})$ and $\rho_{e} = \mathcal{T}(\phi_{e})$ occurring at each re-design step. For the sensitivity filter, the growth is limited to a distance of $r_{\rm min}$ from the original topological boundary, since the information propagates with the operation $d c/ d \phi_{e} = \mathcal{T}^{-1}(d c/d \rho_{e})$ only. On the other hand, the potential growth reaches a distance of $2 r_{\rm min}$ from the topological boundary when using linear density filtering, threshold projection or \ac{HPM}, since the operation $d c/ d \phi_{e} = \mathcal{T}^{-1}(d c/d \rho_{e})$ propagates the sensitivity information to the distance $r_{\rm min}$, and this is further extended by the operation $\rho_{e} = \mathcal{T}(\phi_{e})$.

It is clear from these plots that for a \ac{SIMP} exponent $\eta > 1$, material is most likely to evolve from the current topological boundaries. Thus, while the approach proposed herein is compatible also with linear density filtering and threshold projection, \ac{HPM} very clearly has a stronger ability to reintroduce material. Indeed, when using the density filter alone the sensitivities have maximum values in the interior of the solid region and dissipate rapidly due to only using the decaying hat weighting function (\autoref{eq:filterHatFunction}). Elements beyond the distance $r_{\rm min}$ have almost zero relative density after the two propagation operations, because of the small values of the weights (see \autoref{fig:SensitivityStudy}(f), row A2). On the other hand, the Heaviside projection amplifies the small sensitivity coefficients that propagates from the topological boundaries, boosting the material reintroduction in the whole reach of $2r_{\rm min}$ (see \autoref{fig:SensitivityStudy}(f), row A4). The effect of the threshold projection depends on the threshold value; here we adopt the popular choice of 0.5. In this case, the projection either raises or annihilates the propagated information, and therefore the extent of the propagated sensitivity (and topological growth) remains localized in a small region near the topological boundary (\autoref{fig:SensitivityStudy}, row A3).

For the case of $\eta = 1$ (non-penalized \ac{SIMP}), the discussion above is no longer valid. In this case $s_{e} = \tilde{s}_{e}$ and thus the contribution of low relative density elements to both the response and to the sensitivities is not necessarily negligible, and removing them from the analysis may hinder design evolution. Also, we note that the plots in \autoref{fig:SensitivityStudy} use the linear weighting for the filter (\autoref{eq:filterHatFunction}), since this is the most commonly used. Alternatively, the use of a uniform weighting function, coupled with the Heaviside projection, would reduce the grey regions at topological boundaries and lead to stronger propagation.

The above discussion clearly demonstrates the logic behind the removal of low relative density elements from the analysis space, and illustrates that these can be naturally and effectively reintroduced when using the \ac{HPM}. To summarize:
\begin{itemize}
 \item Whenever $\eta > 2$ is used in the \ac{SIMP} penalization, the contribution of low relative density regions (e.g., $\rho_{e} < 0.1$) to both the structural response and its derivative becomes negligible. Thus, these elements can be removed from the analysis.
 \item The removed elements can be reintroduced in a subsequent optimization step, and this reintroduction is purely governed by the formal sensitivity propagation effect due to the filter, within a distance of $r_{\rm min}$ (sensitivity filter), or $2r_{\rm min}$ (density filter, threshold projection, \ac{HPM}) from the topological boundary.
 \item The material reintroduction occurs only around the topological boundaries, which is anyway the case for a \ac{SIMP} exponent $\eta > 2$. The \ac{HPM}, boosting the sensitivity values along the boundaries, allows a much faster reintroduction of materials into the removed elements than both the linear density filter alone and threshold projection.
\end{itemize}
%

\section{Adaptive nodal boundary conditions}
\label{sec:NodalBoundaryCondition}
The proposed adaptive nodal boundary conditions scheme is applied to a fixed finite element discretization as in past works \cite{BT:03, Lotfi:14}. The basic premise in these works, which optimize mechanical properties, is that if an element density is below a threshold magnitude, here denoted by \enquote{$\rho_{t}$} it can be \enquote{removed} from the state analysis. This \enquote{removal} has taken various forms, including approximating artificial boundary conditions numerically \cite{BT:03} or inserting them explicitly \cite{Lotfi:14}, or neglecting the behavior of \ac{DOF} associated with the element in the nonlinear analysis convergence criteria \cite{PBS:01}. The insertion of boundary conditions effectively eliminates strain energies from appearing in the removed elements. Thus the derivatives of stiffness-related metrics with respect to removed elemental densities are zero. However, the derivative of these metrics with respect to the independent design variables may be nonzero due to the filtering approach as shown in \autoref{fig:SensitivityStudy}(f).

\begin{figure*}[tb]
\centering
  \subfloat[(a)]{
  \includegraphics[width=0.40\linewidth]{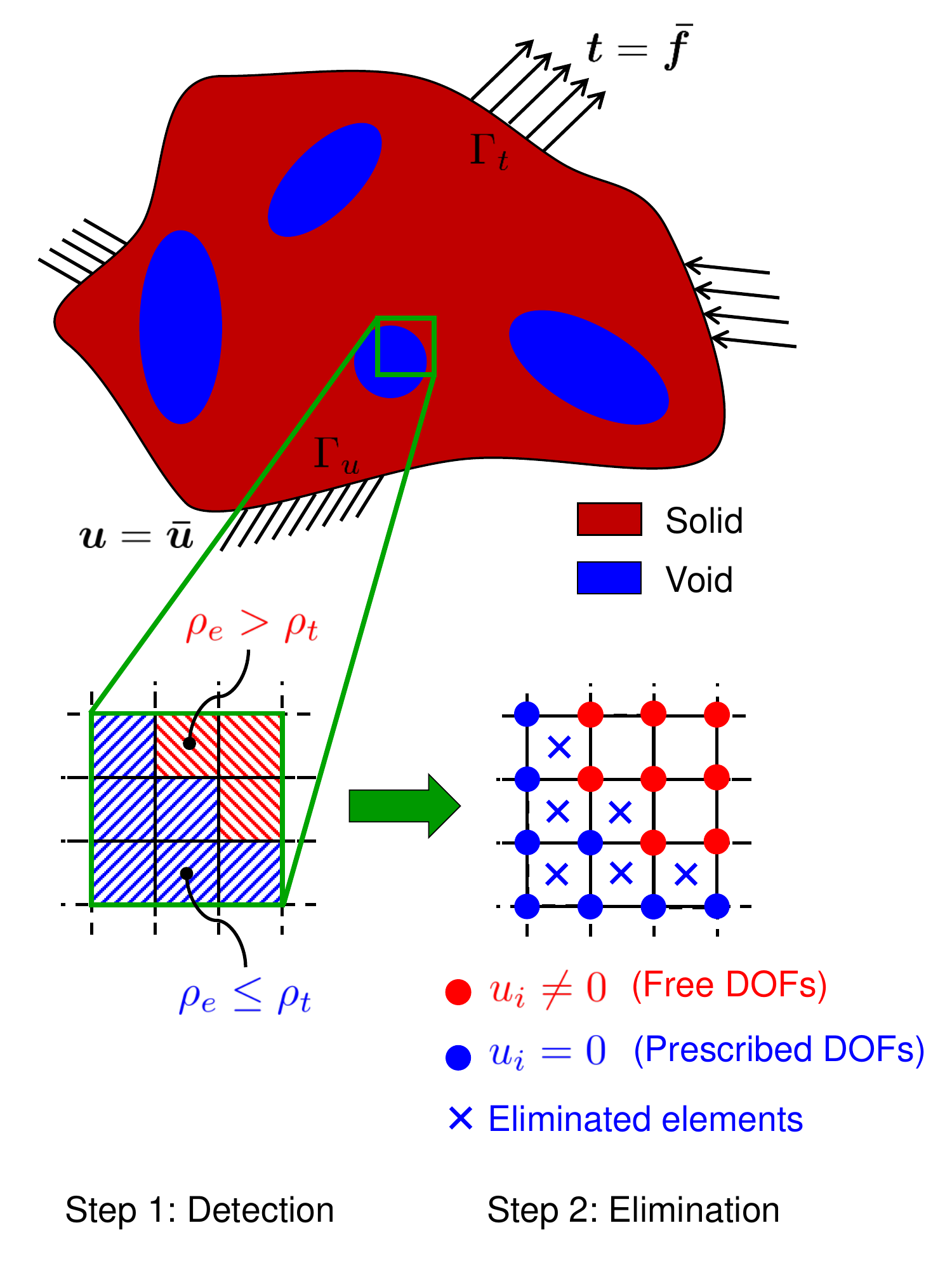}}
  \qquad
  \subfloat[(b)]{
  \includegraphics[width=0.50\linewidth]{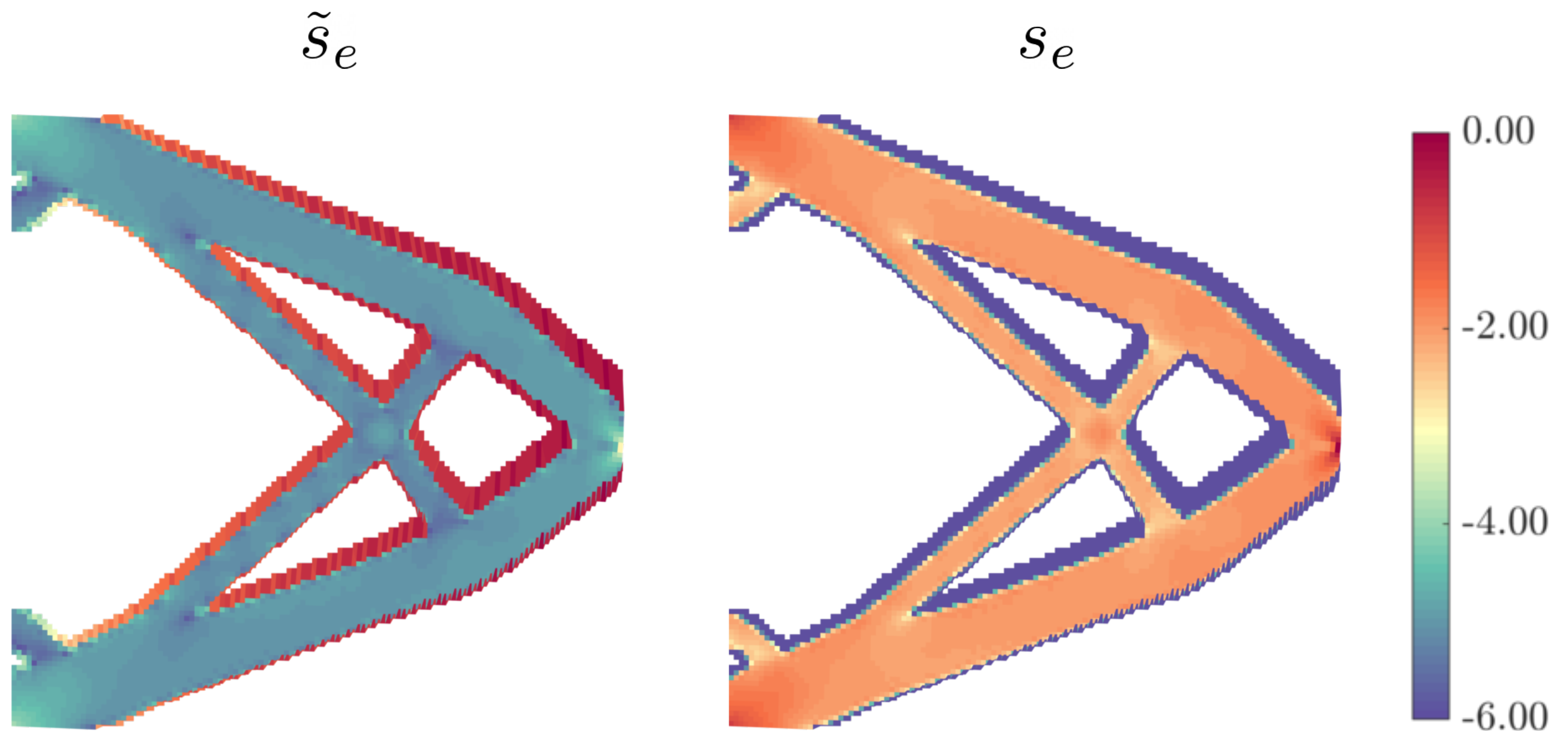}}
\caption{(a) Schematic representation of the element removal approach. In the detection step, the elements where the relative density is below the threshold value $\rho_{t}$ are identified and nodes completely surrounded by these elements are recovered through the connectivity operator. Artificial boundary conditions are then applied to these nodes in the elimination step. (b) The (scaled) deformed configuration of a short cantilever when the nodal boundary conditions are applied. The color map shows the logarithmic scale of the non-dimensional elemental strain energy $\tilde{s}_{e}$ and the dimensional one $s_{e}$, both normalized with respect to their maximum values.}
\label{fig:NoSpliNodalBoundaryScheme}
\end{figure*}

Herein we apply this concept to implement adaptive nodal boundary conditions on displacements in the void regions in structural problems. The approach consists of detection and elimination steps. First, the elements that are deemed to be \enquote{void} are identified, as determined by a threshold magnitude ($\rho_{e} \leq \rho_{t}$). The elimination step then eliminates \ac{DOFs} associated with those elements from the finite element model. As shown in \autoref{fig:NoSpliNodalBoundaryScheme}, this includes all displacement \ac{DOFs} defined at nodes connecting only identified void elements. Displacement \acs{DOFs} on the void-solid interface are left as free \ac{DOFs} (see \autoref{fig:NoSpliNodalBoundaryScheme}(b)). Eliminating \ac{DOFs} in this manner is a straightforward process and reduces the size of the equilibrium system to be solved, leading to computational savings. Additionally, the removed finite elements need not be assembled into the global stiffness matrix and elemental state matrices, such as strain energy, should not be computed. 

A key feature of the proposed approach is that void elements that have been removed may be introduced back into their elemental domain through projection and may readily return as part of the structural topology and their \ac{DOFs} to the analysis model. The numerical examples presented herein will confirm that this occurs. Fundamentally, this property is due to the fact that the design and analysis spaces are separate and discretized on different meshes \cite{GS:10}, and related through projection that occurs over a physical length scale. A key point of emphasis is that this material reintroduction is driven purely by formal sensitivity analysis, as clearly illustrated in Section \ref{sec:SensitivityStudy}.  
 
Of course, it must be acknowledged that formally implementing nodal boundary conditions as the design evolves introduces a discontinuity in the otherwise continuous topology optimization problem, and that designs are likely to be dependent on the chosen magnitude of the threshold $\rho_{t}$. In the case of solid mechanics problems, different threshold magnitudes may lead to different solutions, although if chosen small enough, these solutions tend to have very similar objective function magnitudes. We note that this approach was recently successfully applied to pressure and velocity boundary conditions in incompressible laminar fluid flow problems \cite{BRG:18}, and we have found the fluid flow problems to be less sensitive to this parameter.

To summarize, the following operations are performed at each design iteration for a given $\phi$ and corresponding $\rho_{e}$. (1) Elements with $\rho_{e} \leq \rho{t}$ are set to $\rho_{e} = 0$ and marked for removal, (2) nodal boundary conditions are added to \ac{DOFs} corresponding only to removed elements, (3) the forward problem is solved, (4) sensitivity analysis is performed, (5) a design change in $\phi$ is computed by \ac{MMA} or \ac{OC}, and (6) $\phi$ is projected onto $\rho_{e}$ to update topology.  This procedure is repeated until convergence, as defined in Section \ref{sec:OptimizationProblem}.
%

%
\section{Numerical examples}
 \label{sec:NumericalExamples}
The element removal and reintroduction strategy detailed in \autoref{sec:NodalBoundaryCondition} is now tested on several structural topology optimization examples.

We start with the simplest problem of linear compliance minimization (\autoref{sec:LinComplMinimization}), and then we extend to the force inverter design problem, which is more prone to material reintroduction (\autoref{sec:ForceInverter}). For both cases, we show the practical insensitivity of the obtained designs over a wide range of threshold density values and the remarkable computational savings.

Next, we move to nonlinear state equations. In \autoref{sec:GeometricNonlinearity}, we consider the minimization of the end compliance for a geometrically nonlinear example, that is the setup where the idea of disregarding low relative density elements was also used \cite{BPS:00}. However, we point out that our approach differs from that of \cite{BPS:00}, where elements with low $\rho_{e}$ were excluded from the convergence criterion only, but were still included in the state analysis. Finally, \autoref{ssSec:maxVibration} and \autoref{ssSec:maxBLF} address topology optimization for the maximization of the natural vibration frequency and buckling the load factor, respectively. We note that these are challenging problems, often showing slow convergence and numerical instabilities, and requiring the solution of a computationally expensive eigenvalue equation at each re-design step.

In the following, we will refer to the case where all the elements within the design domain are modeled in the forward analysis as the \enquote{standard approach}. The corresponding design is called the \enquote{reference design}, and the corresponding threshold is denoted by $\rho_{t} = n.a.$ (not applicable). We remark that in this case we will consider a non-zero lower bound $E_{\rm min}$ for the stiffness interpolation, whereas we will set $E_{\rm min} = 0$ when adopting a threshold value $\rho_{t} \geq 0$ since the inserted boundary conditions prevent singularities in the global stiffness matrix.

In all the following numerical examples, $V_{\rm max}$ represents the maximum allowed value of the volume fraction. Since we are using uniform discretizations, the volume of all the elements is the same ($|\Omega_{e}| = h^{2}_{e}$ in 2D and $|\Omega_{e}| = h^{3}_{e}$ in 3D) and the volume fraction is computed simply as
\begin{equation}
 \label{eq:volumeFraction}
  V(\boldsymbol{\phi}) = 
  \frac{1}{m}\sum^{m}_{e =1} \rho_{e}(\boldsymbol{\phi}),
\end{equation}

Finally, in this work we have adopted direct solvers for solving the linear and eigenvalue equations in the forward and adjoint analyses. However, the element removal strategy described in \autoref{sec:NodalBoundaryCondition} is directly compatible with iterative solvers \cite{LHZ+:18}. The following linear elastic designs have been run on a laptop equipped with an Intel(R) Core(TM) i7-5500U@2.40GHz CPU, 16GB of RAM and Matlab 2018b, running in serial mode.

\begin{figure*}[tb]
 \centering
  \subfloat[]{
   \includegraphics[width=0.475\linewidth]
    {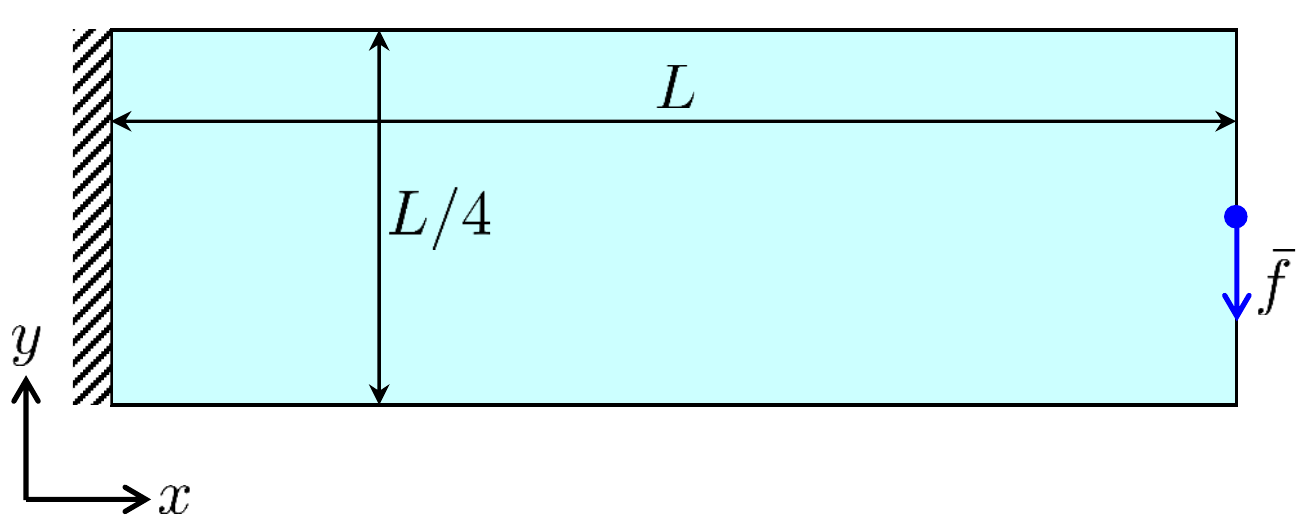}
   }
  \qquad
  \subfloat[]{
   \includegraphics[width=0.40\linewidth]
    {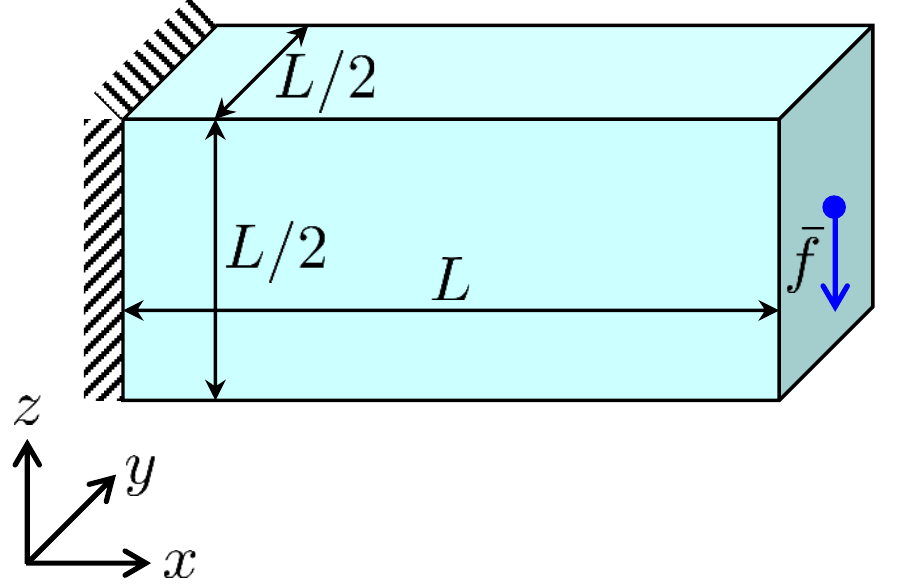}
   }
  \caption{Geometrical setup for the cantilever examples used for the linear and nonlinear compliance minimization examples of \autoref{sec:LinComplMinimization} and \autoref{sec:GeometricNonlinearity}.}
 \label{fig:Cantilever_Beam_LinElast_Schematic}
\end{figure*}

\subsection{Linear elastic design}
\subsubsection{Minimization of the linearized compliance}
 \label{sec:LinComplMinimization}

We consider both 2D and 3D cantilever beams sketched in \autoref{fig:Cantilever_Beam_LinElast_Schematic}, loaded by a concentrated tip force with magnitude $|\boldsymbol{\bar{f}}| = 10^{7}$. We adopt an isotropic material with Poisson's ratio $\nu = 0.3$ and the two-material \ac{SIMP} is used for the Young's modulus interpolation
\begin{equation}
 \label{eq:linCompl_SIMP}
  E(\rho_{e}) = E_{\rm min} + (E_{0}-E_{\rm min})\rho^{\eta}_{e},
\end{equation}
where $E_{0} = 200 \times 10^{9}$ is the modulus of the solid phase and $E_{\rm min} = 10^{-6} E_{0}$ is that of the void phase for the standard approach.

\autoref{eq:OptimizationProblem_1} particularizes to
\begin{equation}
 \label{eq:OptimizationProblem_LinCompliance}
  \begin{aligned}
   \min_{\boldsymbol{\phi}}
   & \: g_{0} = \boldsymbol{u}^{T}\boldsymbol{f} \\
   \text{s.t.}
   & \: K(\boldsymbol{\phi})\boldsymbol{u} - \boldsymbol{f} = \boldsymbol{0} \\
   & V(\boldsymbol{\phi}) - V_{\rm max} \leq 0 \\
   & 0 \leq \phi_{i} \leq 1 \qquad i = 1, \ldots, m   
  \end{aligned},
\end{equation}
where $\partial_{\boldsymbol{u}} \boldsymbol{r}(\boldsymbol{\phi},\boldsymbol{u}(\boldsymbol{\phi})) = K(\boldsymbol{\phi})\boldsymbol{u} - \boldsymbol{f} = \boldsymbol{0}$ are the linear equilibrium equations solved for the (small) displacements $\boldsymbol{u}$. Having a self-adjoint problem, we conclude $\boldsymbol{a} = -\boldsymbol{u}$ in \autoref{eq:DO_SensitivityAnalysis_2} and the sensitivity of the objective becomes $\partial_{\boldsymbol{\rho}}g_{0} = -\boldsymbol{u}^{T}\partial_{\boldsymbol{\rho}}K(\boldsymbol{\phi})\boldsymbol{u}$.

The design update is performed by means of the Optimality Criterion implemented as in \cite{FS:20} and the optimization is considered converged when
\begin{equation}
 \label{eq:D0_terminationCriteria}
 \begin{aligned}
  \tau_{\phi} & = 
  \|\boldsymbol{\phi}^{(k)}-\boldsymbol{\phi}^{(k-1)}\|_{2} / 
  \sqrt{m} < 10^{-4} \\
  \tau_{0} & = (g^{(k)}_{0}-g^{(k-1)}_{0})/g^{(1)}_{0} < 10^{-8}
  \end{aligned},
\end{equation}
hold simultaneously.

For the 2D modeling, we assume plane stress and a fairly large discretization consisting of 1,600 $\times$ 400 $\mathcal{Q}_{4}$ isoparametric elements, for about $1.284 \times 10^{6}$ total \ac{DOFs}. The filter radius, curvature regularization factor and the allowable maximum volume fraction are set to $r_{\rm min} = 7.5 \times 10^{-3}L$ (12 elements), $\beta = 32$ and $V_{\rm max} = 0.5$, respectively. The \ac{SIMP} penalization is started at $\eta = 2$ and increased every 25 optimization steps by $\Delta\eta = 0.5$, up to $\eta_{\rm max} = 6$. For the 3D configuration (see \autoref{fig:Cantilever_Beam_LinElast_Schematic}(b)) we adopt a discretization with $120\times 40 \times 40$ tri-linear $\mathcal{H}_{8}$ elements, for a total of about 610,200 \ac{DOFs}. We set $r_{\rm min} = 2.9 \times 10^{-2}L$ ($2\sqrt{3}$ elements), $\beta = 32$ and the \ac{SIMP} penalization factor is set to $\eta = 3$. For the 3D example, We consider three volume fractions namely $V_{\rm max} = 0.16$, $0.08$ and $0.04$.

\begin{figure*}[tb]
 \centering
  \subfloat[(a) Reference design ($\rho_{t} = n.a.$)]{
   \includegraphics[width = 0.275\linewidth]
    {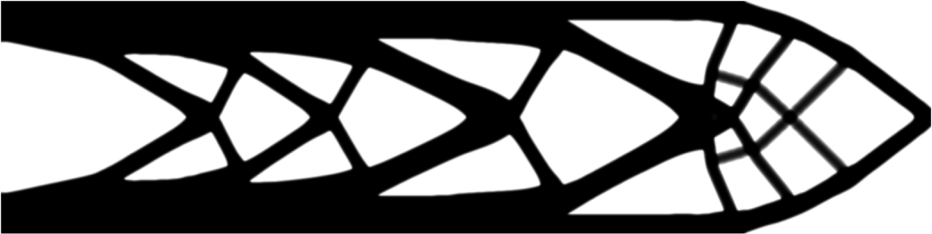}
    }
  \quad
  \subfloat[(b) $\rho_{t} = 0.01$]{
   \includegraphics[width = 0.275\linewidth]
   {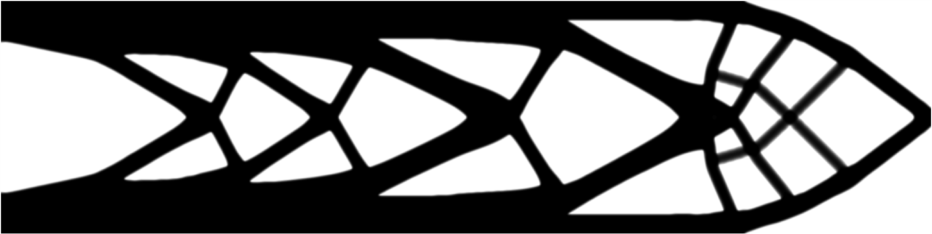}
   }
  \quad
  \subfloat[(c) $\rho_{t} = 0.1$]{
   \includegraphics[width = 0.275\linewidth]
   {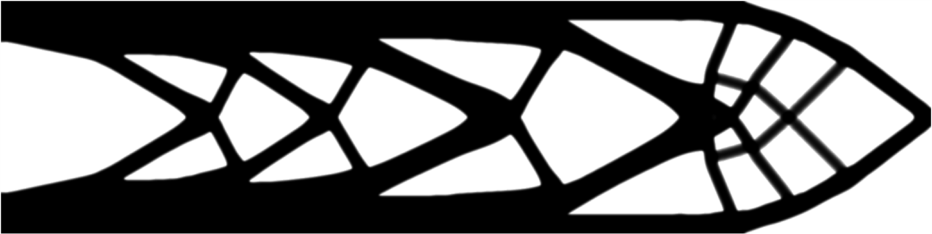}
   }
  \\
  \subfloat[(d) $V_{\rm max} = 0.16$]{
   \includegraphics[width = 0.275\linewidth]
   {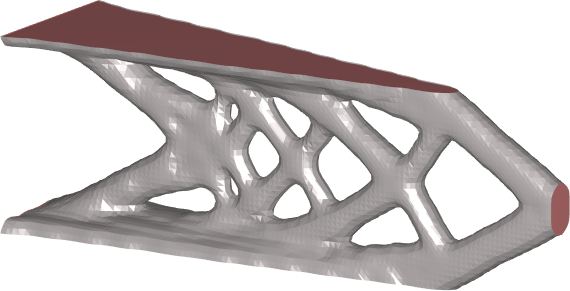}
   }
   \quad
   \subfloat[(e) $V_{\rm max} = 0.08$]{
   \includegraphics[width = 0.275\linewidth]
   {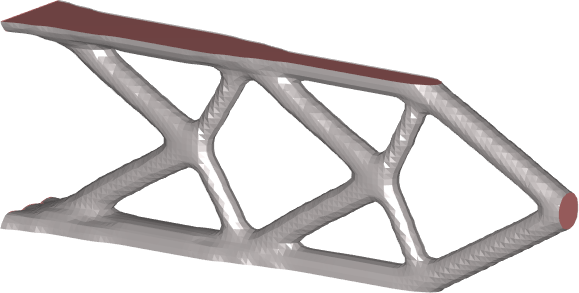}
   }
   \quad
  \subfloat[(f) $V_{\rm max} = 0.04$]{
   \includegraphics[width = 0.275\linewidth]
   {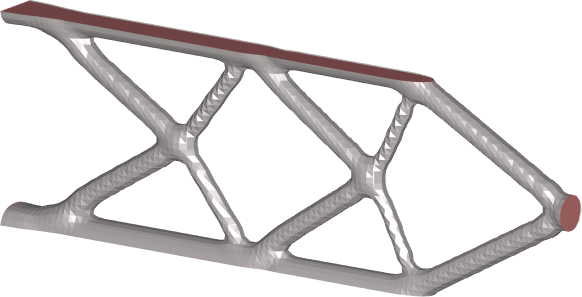}
   }
  \caption{Linear compliance minimization: 2D designs corresponding to different threshold densities (a-c), and 3D designs with different volume fractions, all corresponding to the threshold density $\rho_{t} = 0.01$ (d-f). The 3D plots have been thresholded at the relative density value $\bar{\rho}_{e} = 0.95$.}
 \label{fig:linCompl_topResults}
\end{figure*}

The reference 2D design is shown in \autoref{fig:linCompl_topResults}(a), whereas plots (b,c) display the designs obtained when using threshold densities $\rho_{t} = 0.01$ and $\rho_{t} = 0.1$, respectively. The designs are basically indistinguishable, and the only small differences happen to be localized at the boundaries, as expected. \autoref{fig:linCompl_topResults}(d)-(f) show the 3D designs corresponding to $\rho_{t} = 0.01$, for the three given volume fractions. These are also hardly distinguishable from the corresponding reference designs and from those obtained with the higher threshold density, i.e., $\rho_{t}= 0.1$.

\begin{figure*}[tb]
 \centering
  \subfloat[(a) 2D]{
   \includegraphics[width = 0.285\linewidth]
    {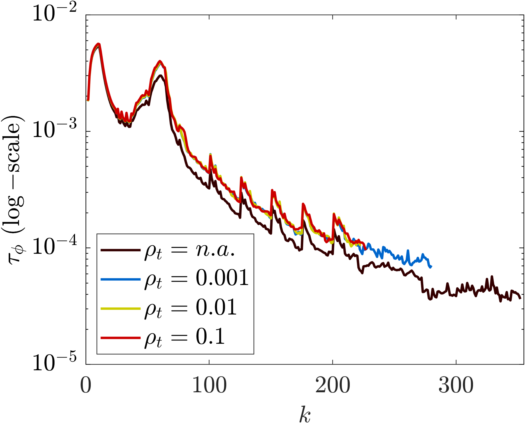}
    }
  \quad
  \subfloat[(b) 2D]{
   \includegraphics[width = 0.305\linewidth]
    {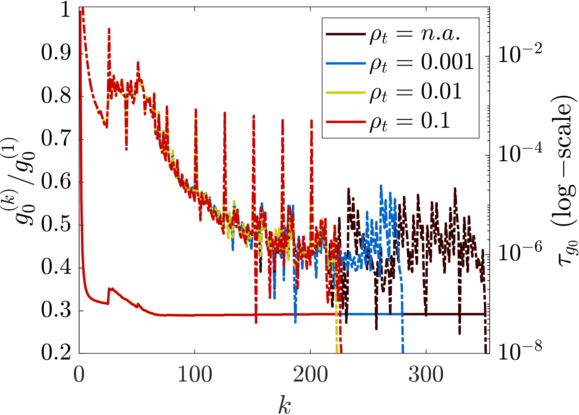}
   }
   \quad
  \subfloat[(c) 2D]{
   \includegraphics[width = 0.305\linewidth]
    {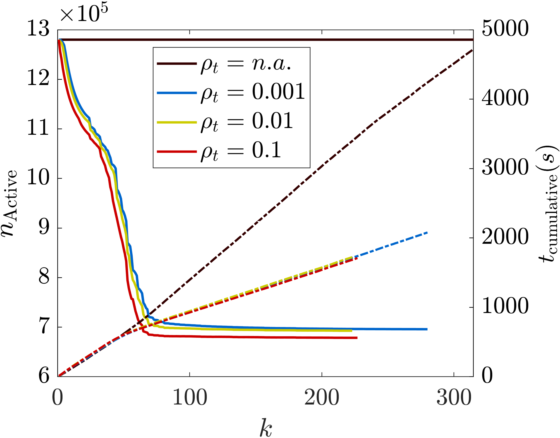}
   }
   \\
     \subfloat[(d) 3D, $V_{\rm max} = 0.16$]{
   \includegraphics[width = 0.285\linewidth]
    {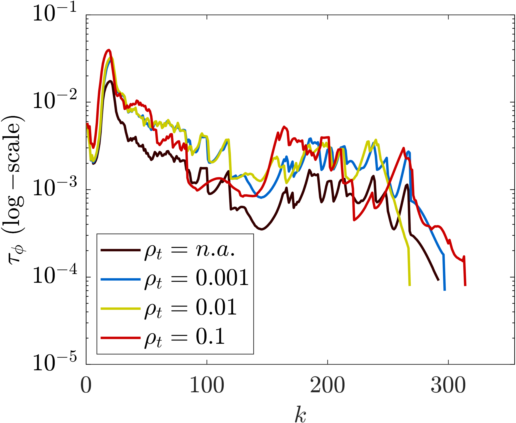}
    }
  \quad
  \subfloat[(e) 3D, $V_{\rm max} = 0.16$]{
   \includegraphics[width = 0.305\linewidth]
    {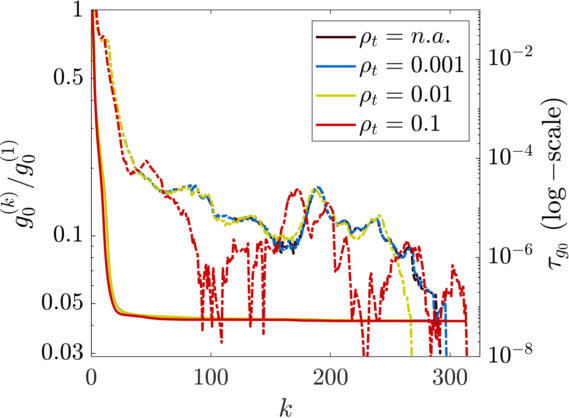}
   }
   \quad
  \subfloat[(f) 3D, $V_{\rm max} = 0.16$]{
   \includegraphics[width = 0.305\linewidth]
    {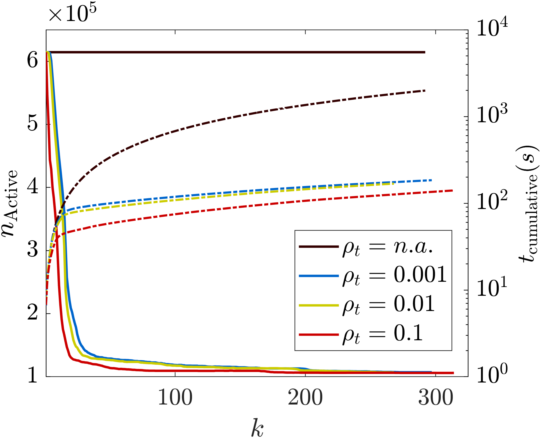}
   }
  \caption{Optimization histories for the 2D and 3D linearized compliance minimization example. Plots (a,d) show the change in the design variables as defined by \autoref{eq:D0_terminationCriteria}, (b,e) show the evolution of the normalized compliance ($g^{(k)}_{0}/g^{(1)}_{0}$, solid lines plotted against the left axis) and the relative change in the objective function as defined by \autoref{eq:D0_terminationCriteria} ($\tau_{g_{0}}$, dashed lines plotted against the right axis). Plots (c,f) display the number of \enquote{active} elements in the state analysis ($n_{\rm Active}$, continuous curves plotted against the left axis) and the cumulative computational time (dashed curves plotted against the right axis). For the 3D example, plot (f) is in $\log$-scale.}
 \label{fig:linCompl_2D_optResults}
\end{figure*}

The evolution of the normalized compliance ($g^{(k)}_{0}/g^{(1)}_{0}$), the parameters $\tau_{g_{0}}$ and $\tau_{\phi}$ measuring the convergence of the design process, and of the number of \enquote{active} \ac{DOFs} ($n_{\rm Active}$) are shown in \autoref{fig:linCompl_2D_optResults}. The reference values for the optimized compliance, normalized with respect to the initial one, are $g^{\ast}_{0,{\rm Ref}} = 0.292026$ and $0.041944$ for the 2D and 3D reference designs, respectively. For the designs corresponding to other thresholds $\rho_{t}$, the relative differences, defined by
\begin{equation}
 \label{eq:relDiffObjFunction}
  \Delta g^{\ast}_{0}(\rho_{t}) =
  \frac{g^{\ast}_{0}(\rho_{t})}
  {g^{\ast}_{0,{\rm Ref}}}-1,
\end{equation}
are listed in \autoref{tab:convergence2D3DlinearizedCompliance}, including other relevant data for the maximum allowed volume fractions $V_{\rm max} = 0.5$ (2D setup) and $V_{\rm max} = 0.16$ (3D setup).

\begin{table*}[!tb]
 \begin{center}
  \caption{Numerical results for the linearized compliance minimization example using different values of the threshold density $\rho_{t}$ with $V_{\rm max} = 0.5$ for 2D and $V_{\rm max} = 0.16$ for 3D. $k_{\ast}$ is the iteration at which the convergence criterion (see \autoref{eq:D0_terminationCriteria}) is fulfilled, $\Delta g_{0}$ is the relative difference in the objective function defined by \autoref{eq:relDiffObjFunction}, $m_{\rm ND}$ is the non-discreteness measure \cite{Sigmund:07} and $t_{\rm avg}$ is the average iteration time.}
   \label{tab:convergence2D3DlinearizedCompliance}
    \begin{tabular}{cccccc}
     \hline\noalign{\smallskip}
     & $\rho_{t}$ & $k_{\ast}$ & $\Delta g^{\ast}_{0}$ & $m_{\rm ND}$ ($10^{-4}$) & $t_{\rm avg} \ \rm (sec)$ \\
     \hline\noalign{\smallskip}
     \multirow{4}{*}{2D} &   $n.a.$ & 352 &         0 &  9.4 & 14.9 \\
                         & $0.001$ & 280 & $10^{-5}$ & 11.5 &  7.7 \\
                         &  $0.01$ & 227 & $10^{-5}$ & 15.5 &  7.5 \\
                         &  $0.1$ & 223 & $10^{-3}$ & 14.6 &  7.4 \\
     \hline\noalign{\smallskip}
     \multirow{4}{*}{3D} &   $n.a.$ & 292 &         0 & 3.8 & 91.2 \\
                         & $0.001$ & 297 & $10^{-5}$ & 3.8 & 10.5 \\
                         &  $0.01$ & 268 & $10^{-4}$ & 1.9 & 10.3 \\
                         &  $0.1$ & 314 & $10^{-3}$ & 3.1 &  5.5 \\
     \noalign{\smallskip}\hline
    \end{tabular}
 \end{center}
\end{table*}

\autoref{fig:linCompl_2D_optResults}(a,b) show that, for the 2D setup, the optimization converges quicker as the threshold density is increased. However, this is not always the case for the 3D setup (see \autoref{fig:linCompl_2D_optResults}(d,e) and also \autoref{tab:convergence2D3DlinearizedCompliance}). The evolution of the compliance corresponding to different threshold densities is essentially indistinguishable from plots \autoref{fig:linCompl_2D_optResults}(b,e). From \autoref{tab:convergence2D3DlinearizedCompliance}, we conclude that the relative difference in the final compliance is on the order of a tenth of a percent or less, for both the 2D and 3D cases.

We refer to \autoref{fig:linCompl_2D_optResults}(c,f) to discuss the computational savings. For the 2D setting, the element removal strategy cuts the number of \enquote{active} \ac{DOFs} to about $85\%$ after 20 re-design steps, $70\%$ after 50 steps and $50\%$ after 100 steps, and this trend is about the same for the three chosen values of $\rho_{t}$. For the 3D setting, referring to the case $V_{\rm max} = 0.16$, for the threshold values $\rho_{t} = [0.001, 0.01, 0.1]$, $n_{\rm Active}$ is cut to $[35\%, 30\%, 23\%]$ after 20 iterations, $[24\%, 23\%, 21\%]$ after 50 iterations and $[23\%, 23\%, 20\%]$ after 100 iterations, respectively. Thus, the reduction in the number of \enquote{active} \ac{DOFs} happens to be quicker than in 2D, and also slightly more sensitive to the threshold density value as can be seen in \autoref{fig:linCompl_2D_optResults}(f).

The cumulative CPU time increases at a nearly linear rate for the standard approach whereas, when the element removal strategy is employed, we notice a remarkable drop as soon as $n_{\rm Active}$ starts decreasing. For the 2D case, the average iteration time of $t_{\rm avg} = 14.9 s$ for the classical approach, reduces to about $t_{\rm avg} = 7.5 s$ for $\rho_{t}\in [0.001, 0.1]$. For the 3D case, the savings are much more substantial, as the average iteration time reduces from $t_{\rm avg} = 91 s$ in the standard approach, to about $t_{\rm avg} = 10.5 s$ for $\rho_{t} = 0.001$, $\rho_{t} = 0.01$ and $t_{\rm avg} = 5.5 s$ for $\rho_{t} = 0.1$ (see \autoref{tab:convergence2D3DlinearizedCompliance}). Savings for the lower volume fractions $V_{\rm max} = 0.08$ and $0.04$ are even larger as more void elements are present in these designs.

\subsubsection{Force inverter}
 \label{sec:ForceInverter}
\begin{figure*}[tb]
 \centering
  \subfloat[(a)]{
   \includegraphics[width=0.55\linewidth]
   {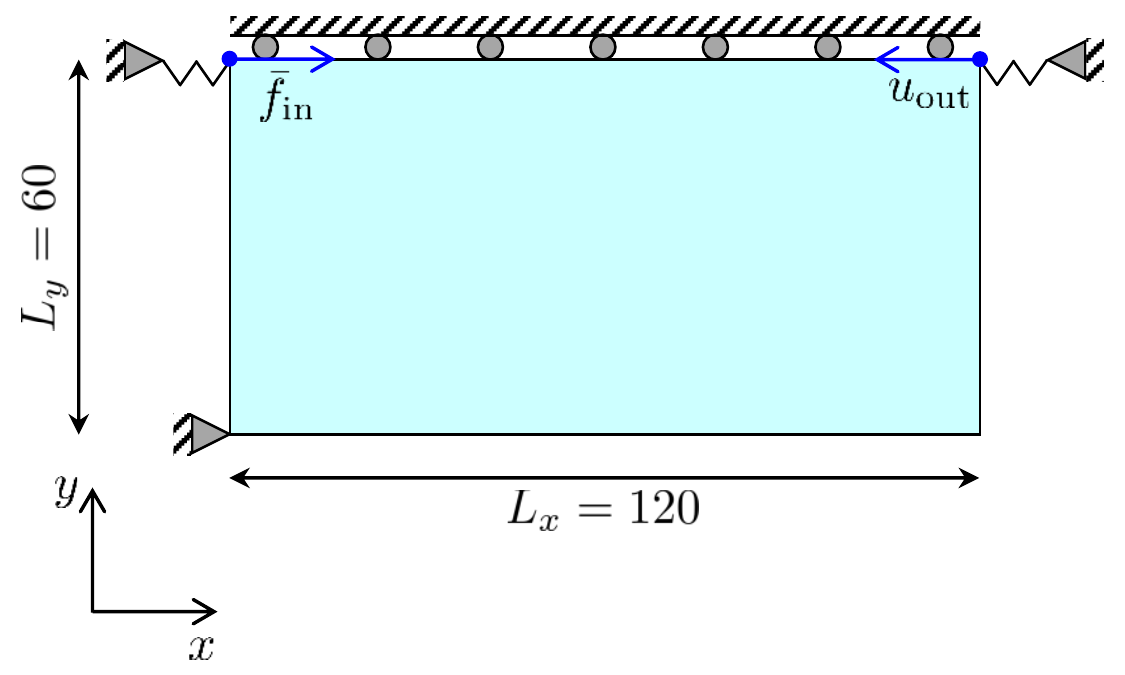}}
   \quad
  \subfloat[(b) $\rho_{t} = n.a.$, $g_{0} = -3.47$]{
   \includegraphics[width=0.30\linewidth]
   {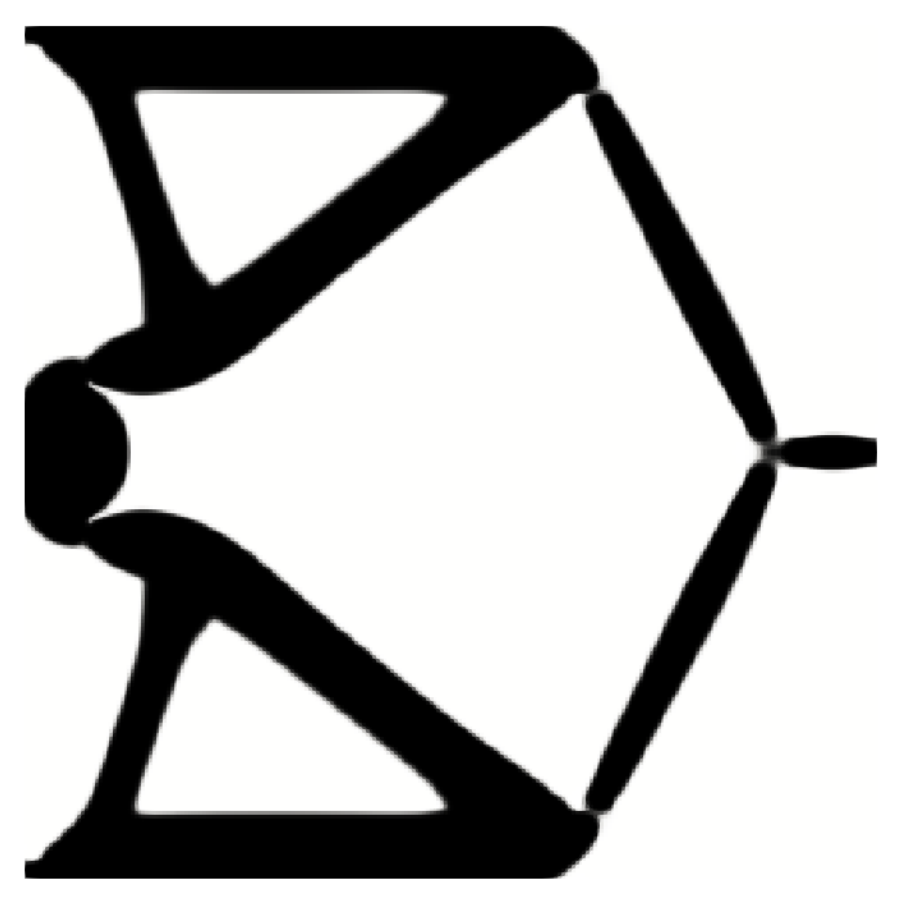}}
  \\
  \subfloat[(c) $\rho_{t} = 0.001$, $g_{0} = -4.86$]{
   \includegraphics[width=0.30\linewidth]
   {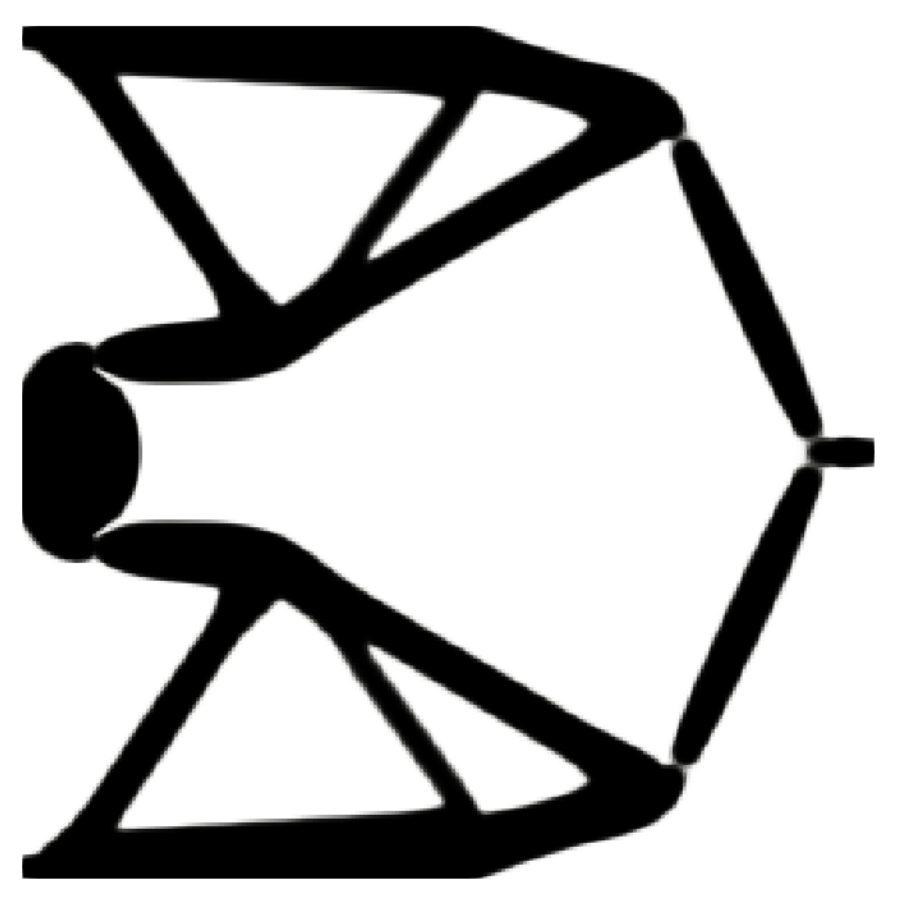}}
   \quad
  \subfloat[(d) $\rho_{t} = 0.01$, $g_{0} = -4.85$]{
   \includegraphics[width=0.30\linewidth]{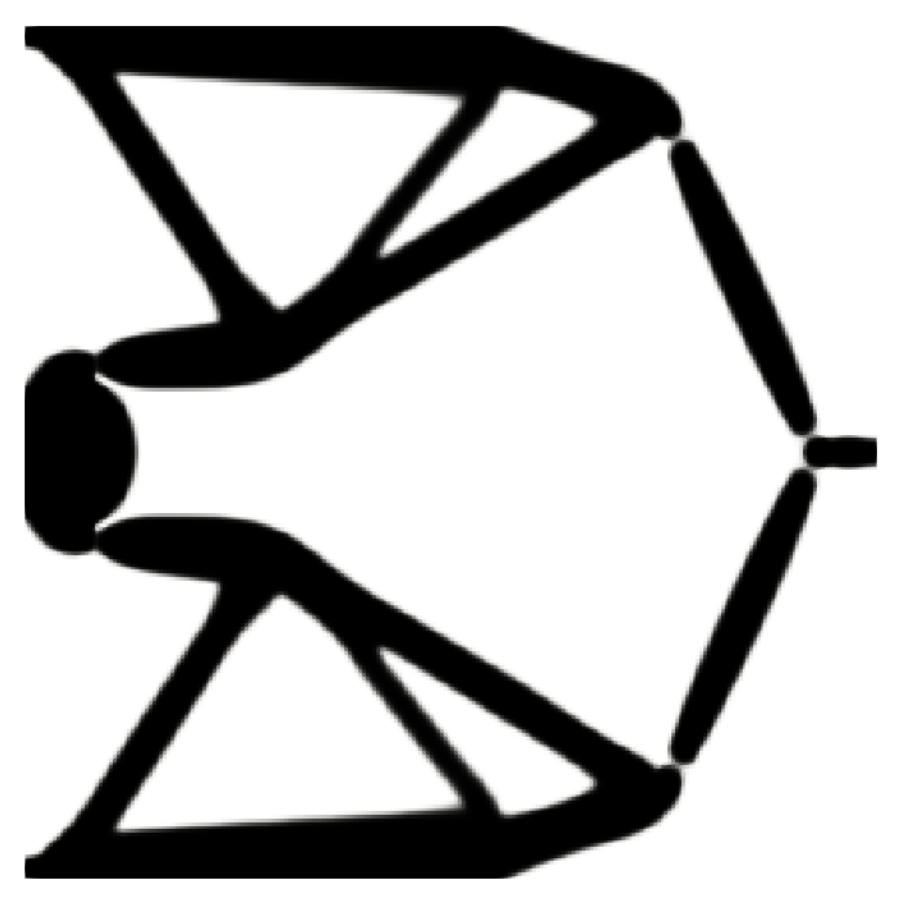}}
  \quad
  \subfloat[(e) $\rho_{t} = 0.1$, $g_{0} = -4.88$]{
   \includegraphics[width=0.30\linewidth]{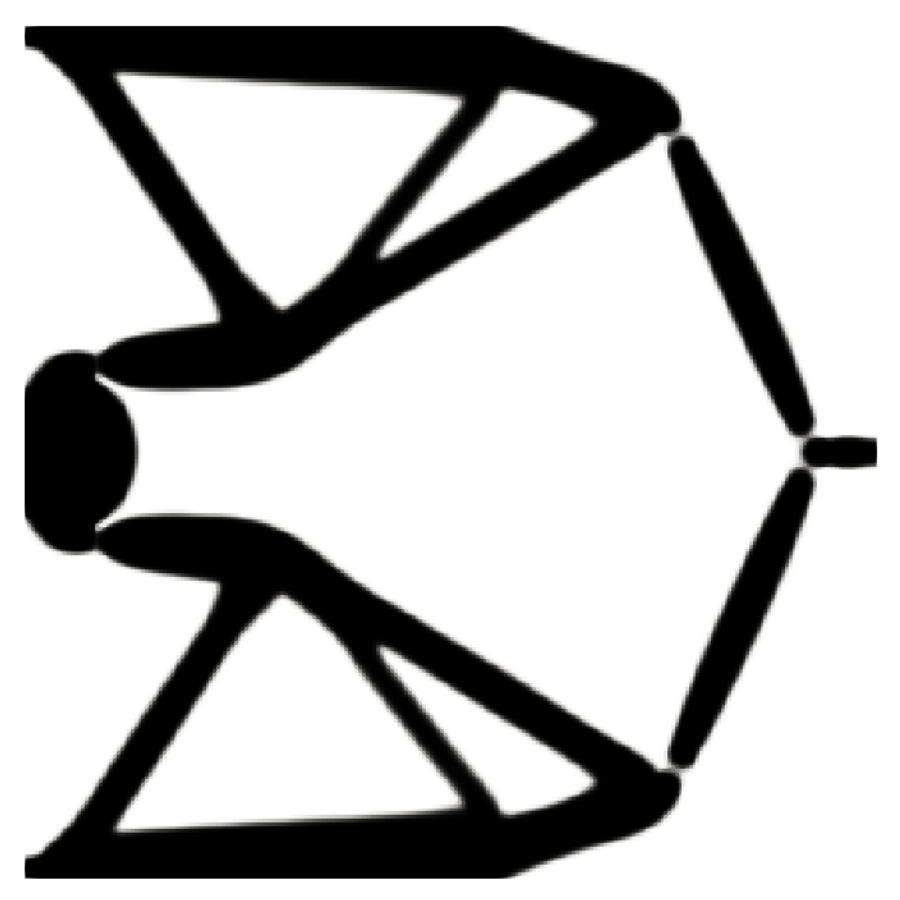}}
\caption{Geometrical setup for the 2D linear elastic force inverter example (only the lower half of the structure is modeled due to symmetry) (a) and optimized designs corresponding to different threshold density values (b-e).}
\label{fig:nbc_n_y_force_inverter_2D}
\end{figure*}

We now consider the design of the force inverter sketched in \autoref{fig:nbc_n_y_force_inverter_2D}(a), where we model only the lower half of the structure due to symmetry around the $x$ axis. Assuming linearly elastic behavior, the objective is to maximize the leftward displacement at the output point ($u_{\text{out}}$) for a given rightward force at the input point ($\bar{f}_{\text{in}}$).

The optimization problem is formulated as follows:
\begin{equation}
 \label{eq:ForceInverterLinearElasticObj}
  \begin{aligned}
   \min_{\boldsymbol{\phi}}
   & \: g_{0} = \boldsymbol{l}^{T}\boldsymbol{u} \\
   \text{s.t.}
   & \: K(\boldsymbol{\phi})\boldsymbol{u} - \boldsymbol{f} = \boldsymbol{0} \\
   & V(\boldsymbol{\phi}) - V_{\rm max} \leq 0 \\
   & 0 \leq \phi_{i} \leq 1 \qquad i = 1, \ldots, N_{i}   
  \end{aligned},
\end{equation}
where $\boldsymbol{l}$ is the vector targeting the output displacement \ac{DOF} (i.e., all $l_{i} = 0$, except for $l_{j} = 1$, where $j$ is the output \ac{DOF}), and $N_{i}$ is the number of design variables.

Compared to compliance minimization, this mechanism design problem poses more challenges to the element removal strategy. Indeed, a trivial solution for the maximization of the negative displacement would be to disconnect the output node from the structure, or to create hinges of soft material. Thus, the evolution of the design towards a well defined structural layout is more complicated, with frequent material removal and reintroduction.

The design domain is discretized with $240 \times 120$ isoparametric $\mathcal{Q}_{4}$ elements. The magnitude of the applied input force is set to $|\bar{f}_{\text{in}}| = 1.0$. Two springs are added to the model at the input and output points with stiffness $k_{\mathrm{in}} = 1.0$ and $k_{\mathrm{out}} = 1.0 \times 10^{-3}$. The Young's modulus and Poisson's ratio of the base materials are set to $E_{0} = 1.0$ and $\nu = 0.3$, respectively, and $E_{\rm min} = 10^{-3} E_{0}$ is that of the void phase when $\rho_{t} = n.a.$. The elemental Young modulus is interpolated using \autoref{eq:linCompl_SIMP}. A continuation on the \ac{SIMP} exponent $\eta$ is performed with tightened \ac{MMA} asymptotes \cite{GAH:11}, starting with $\eta = 2$ and increasing the penalization by $\Delta \eta  = 0.2$ at every $20$ design iterations up to $\eta_{max} = 8$. The filter radius, the curvature of regularization, and the prescribed volume fraction are set to $r_{\rm min} = 1$, $\beta = 50$ and $V_{\rm max} = 0.3$, respectively.

The reference design ($\rho_{t} = n.a.$), obtained after $800$ optimization steps, is shown in \autoref{fig:nbc_n_y_force_inverter_2D}(b). The initial value of the objective function is $g^{(1)}_{0}(\rho_{t} = n.a.) = 0.4343$, and for the reference design we get $g^{\ast}_{0}(\rho_{t} = n.a.) = -3.47$.

\autoref{fig:nbc_n_y_force_inverter_2D}(c-e) show the designs obtained when using threshold densities $\rho_{t} = 0.001$, $\rho_{t} = 0.01$, and $\rho_{t}=0.1$, respectively. As the threshold density increases, the optimized design show a slight deviation from the reference design. We notice that the absolute value of the output displacement is considerably higher when using the element removal approach. This happens because in the reference design (\autoref{fig:nbc_n_y_force_inverter_2D}(b)) the \enquote{void} elements in the inner part of the domain are treated numerically as soft elements with Young's modulus of $E_{\rm min}$, which are restraining, to some extent, the activation of the mechanism connecting the output node. Removing those elements clearly allows the output node to undergo a larger horizontal displacement. This is a good example of the sensitivity of mechanism design problems.

\begin{figure*}[!tb]
 \centering
  \subfloat[(a)]{
   \includegraphics[width = 0.435\linewidth]
    {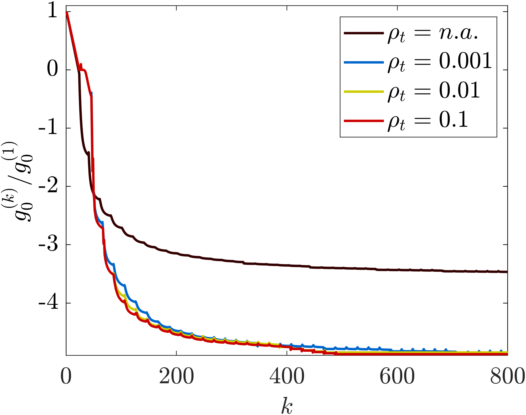}
   }
   \qquad
  \subfloat[(b)]{
   \includegraphics[width = 0.45\linewidth]
    {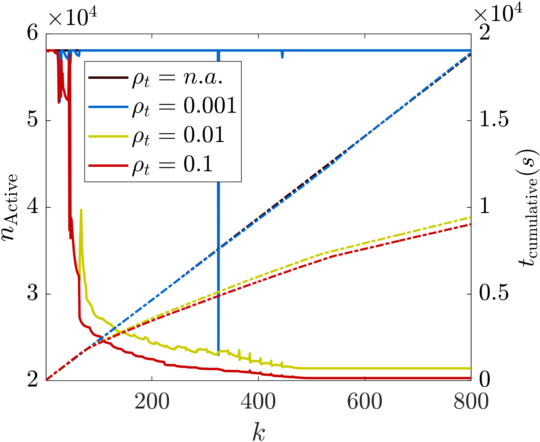}
   }
  \caption{Optimization history for the force inverter problem: (a) evolution of the normalized compliance ($g^{(k)}_{0}/g^{(1)}_{0}$) and (b) the number of \enquote{active} \ac{DOFs} in the state analysis ( continuous curves plotted against the left axis) and the cumulative computational time (dashed curves plotted against the right axis). An increase in the number of \enquote{active} \ac{DOFs} indicates previously removed elements have been reintroduced by the optimizer.}
 \label{fig:forceInverter_2D_optResults}
\end{figure*}

\autoref{fig:forceInverter_2D_optResults} shows the evolution of the objective function, normalized with respect to its initial value, for the different threshold densities. The jumps in the design history correspond to the update in the \ac{SIMP} exponent during the continuation procedure. We observe a marked difference in the evolution of the objective functions. For the threshold density $\rho_{t} > 0$, the objective is always below the one from the standard approach. This clearly indicates numerical artifact of low-density elements on the structural response that cause significant changes in the objective function, for $\rho_{t} > 0$.

It is clear from \autoref{fig:forceInverter_2D_optResults}(b) that the number of \enquote{active} \ac{DOFs} and progressive CPU time are both reduced as the optimization progresses. It is also clear that material is readily reintroduced, evidenced by the jump in free \ac{DOFs}. The use of the low threshold value $\rho_{t} = 0.01$ brings no savings, because the element reintroduction is very frequent and all the elements are active for the most of the optimization history. For larger $\rho_{t}$, on the other hand, the computational savings become apparent.

\subsection{Geometric nonlinearity
\label{sec:GeometricNonlinearity}}
In this section, topology optimization of structures undergoing large deformations is considered by solving the 2D benchmark optimization problem studied in \cite{BPS:00}. The schematic representation of the 2D problem is given in \autoref{fig:Cantilever_Beam_LinElast_Schematic} where the length and thickness are set to $L = 1$ and $0.1$, respectively. Different magnitudes for the applied point load $\bar{f}$ are considered, i.e., $\bar{f} = [12, 144, 240] \times 10^{3}$ as in \cite{BPS:00}. 

The design domain is discretized with $400 \times 100$ isoparametric $\mathcal{Q}_{4}$ elements. The elemental stiffness of the structure is interpolated using \autoref{eq:linCompl_SIMP}, where $E_{\text{min}} = 10^{-6} E_{0}$ for $\rho_{t} = n.a.$. The Young's modulus and Poisson's ratio of the base materials are set to $E_{0} = 3 \times 10^{9}$ and $\nu = 0.4$, respectively, following \cite{BPS:00}. A continuation on the \ac{SIMP} exponent $\eta$ is performed, starting with $\eta = 2.5$ and increasing the penalization by $\Delta \eta  = 0.5$. We note that we consider $400$ design iterations for the first continuation step to encourage the optimization problem to converge to a stable solution. For the rest of continuation steps, $\Delta \eta$ is increased at every $50$ design iterations. The filter radius  and the curvature of regularization are set to $r_{\rm min} = 0.008$ and $\beta = 50$, respectively. The prescribed volume fraction for the solid materials is set to $V_{\rm max} = 0.5$.

The geometric nonlinear behavior of the structure is modeled using the total Lagrangian finite element formulation \cite{BLM+:13}. Using the principle of virtual work, the weak form of the governing equation, in the absence of body forces and accelerations, can be written as:
\begin{equation}\label{eq:StructuralNLinearElasticObj_2}
\begin{split}
\boldsymbol{r} (\boldsymbol{\phi}, \boldsymbol{u} (\boldsymbol{\phi})) = \int_{\Omega_{0}} \boldsymbol{S} \colon \delta  \boldsymbol{E} \ d\Omega - \int_{\Gamma_{0}} \boldsymbol{\bar T} \cdot \delta \boldsymbol{u} \ d\Gamma = 0,
\end{split}
\end{equation}
where $\boldsymbol{S}$ is the second Piola-Kirchhoff stress tensor, $\boldsymbol{E}$ is the nonlinear Green-Lagrange strain tensor, $\delta \boldsymbol{u}$ is the virtual displacement field, and $\boldsymbol{\bar T}$ is the traction in the reference configuration \cite{BLM+:13}. The nonlinear Green-Lagrange strain tensor is defined as follows:
\begin{equation}\label{eq:GeomNlnr_GreenStrain}
\begin{split}
\boldsymbol{E} = \frac{1}{2} \big( \boldsymbol{F}^{T} \boldsymbol{F} - \boldsymbol{I} \big),
\end{split}
\end{equation}
where $\boldsymbol{F} = \nabla_{0} \boldsymbol{u} + \boldsymbol{I}$ is the deformation gradient tensor with $\nabla_{0}$ that represents the  gradient  with  respect  to  the  reference coordinates and $\boldsymbol{I}$ is the identity tensor.

The second Piola–Kirchhoff stress tensor is obtained from the modified Saint Venant strain energy function, which is shown suitable for materials undergoing excessive compression \cite{KS:13}. The second Piola–Kirchhoff stress tensor and the strain energy function are given as follows:
\begin{equation}\label{eq:GeomNlnr_PK2_SaintVenant}
\begin{split}
\boldsymbol{S} = 2 \frac{\partial W}{\partial \boldsymbol{C}}, \quad \text{with} \quad W = \frac{1}{2} \bar{\lambda} \big( J - 1 \big)^{2} + \bar{\mu} \text{tr} (\boldsymbol{E}^{2}),
\end{split}
\end{equation}
where $\boldsymbol{C} = \boldsymbol{F}^{T} \boldsymbol{F}$ is the right Cauchy–Green deformation tensor, $J = \det \boldsymbol{F}$, and and $\bar{\lambda}$ and $\bar{\mu}$ are Lame's material parameters.

The optimization problem is formulated as follows:
\begin{equation}
 \label{eq:GeoNlnrOptProblem}
  \begin{aligned}
   \min_{\boldsymbol{\phi}}
   & \: g_{0} = \boldsymbol{F}_{\text{ext}}^{T} \boldsymbol{\hat{u}} \\
   \text{s.t.}
   & \: \boldsymbol{r} (\boldsymbol{\phi}, \boldsymbol{u} (\boldsymbol{\phi})) = \boldsymbol{0} \\
   & V(\boldsymbol{\phi}) - V_{\rm max} \leq 0 \\
   & 0 \leq \phi_{i} \leq 1 \qquad i = 1, \ldots, N_{i}   
  \end{aligned},
\end{equation}
where $\boldsymbol{F}_{\text{ext}}$ is the external force vector, $\boldsymbol{\hat{u}}$ is the displacement solution of \autoref{eq:StructuralNLinearElasticObj_2} in the reference configuration, and $N_{i}$ is the number of design variables.

\begin{figure*}[tb]
    \setlength{\tabcolsep}{1pt}
    \begin{tabular}{ccc}
    \centering
    \subfloat[$\bar{f} = 12$ kN, $\rho_{t} = n.a.$]{
    \includegraphics[width=0.32\linewidth]{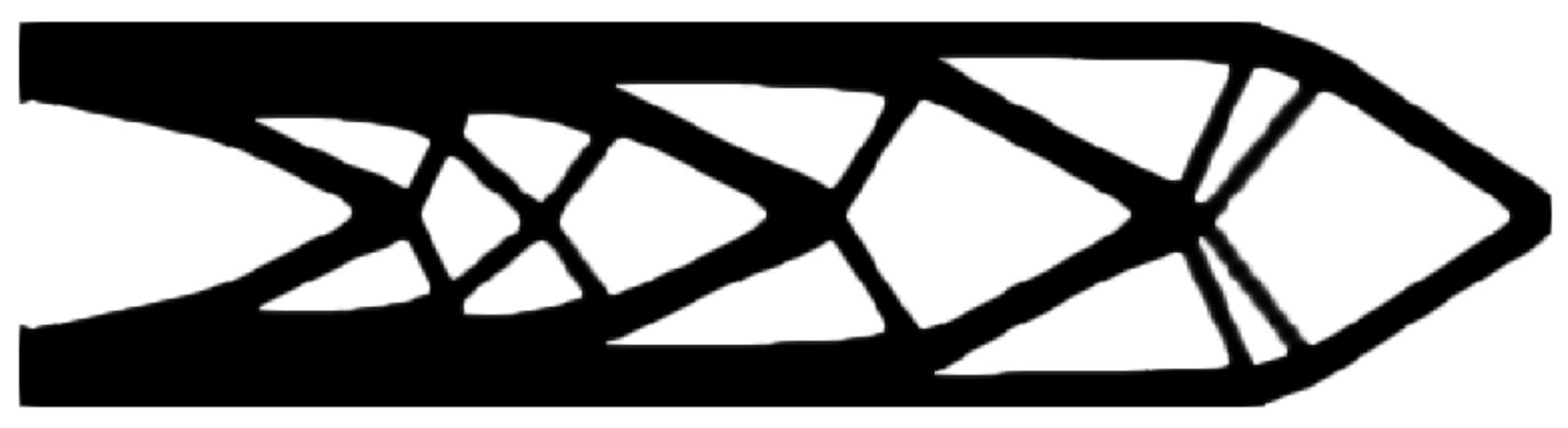} } &
    \subfloat[$\bar{f} = 144$ kN, $\rho_{t} = n.a.$]{
    \includegraphics[width=0.32\linewidth]{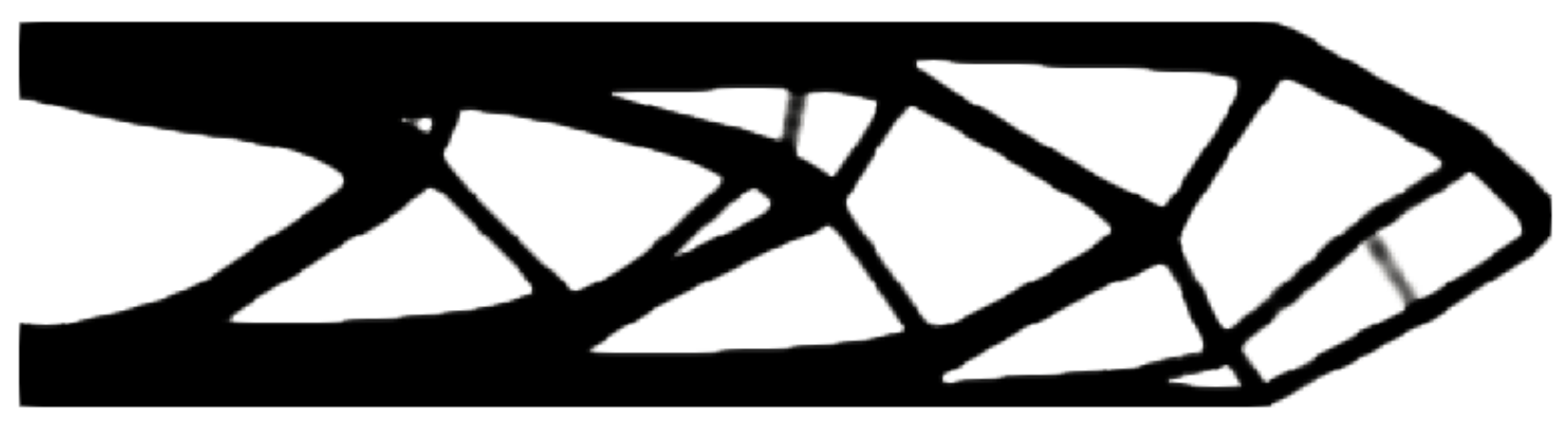} } & \\
    \newline
    \subfloat[$\bar{f} = 12$ kN, $\rho_{t} = 0.01$]{
    \includegraphics[width=0.32\linewidth]{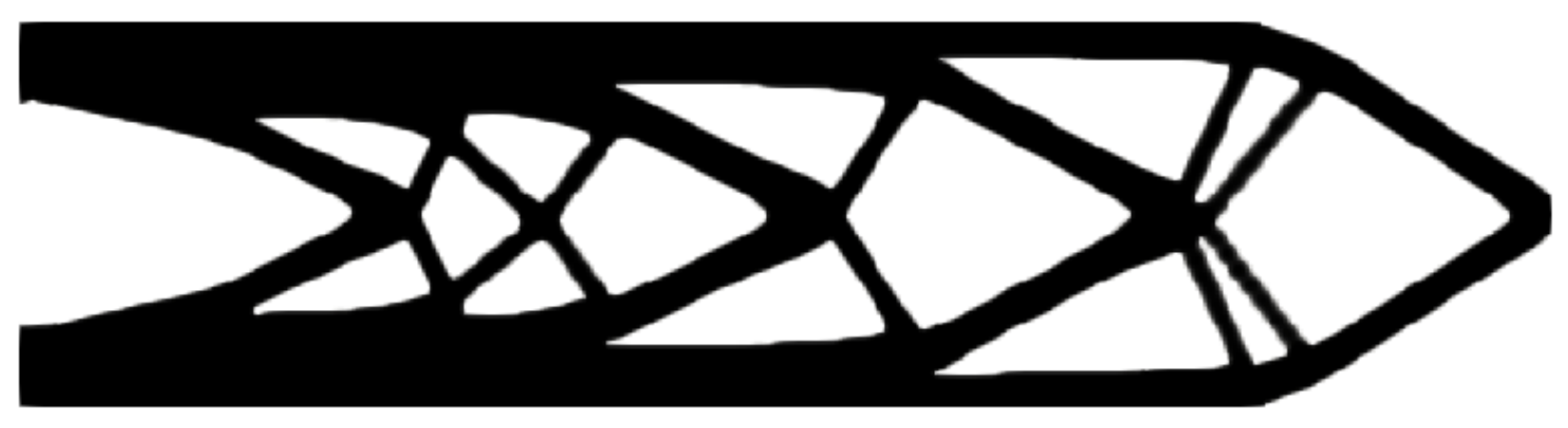} } &
    \subfloat[$\bar{f} = 144$ kN, $\rho_{t} = 0.01$]{
    \includegraphics[width=0.32\linewidth]{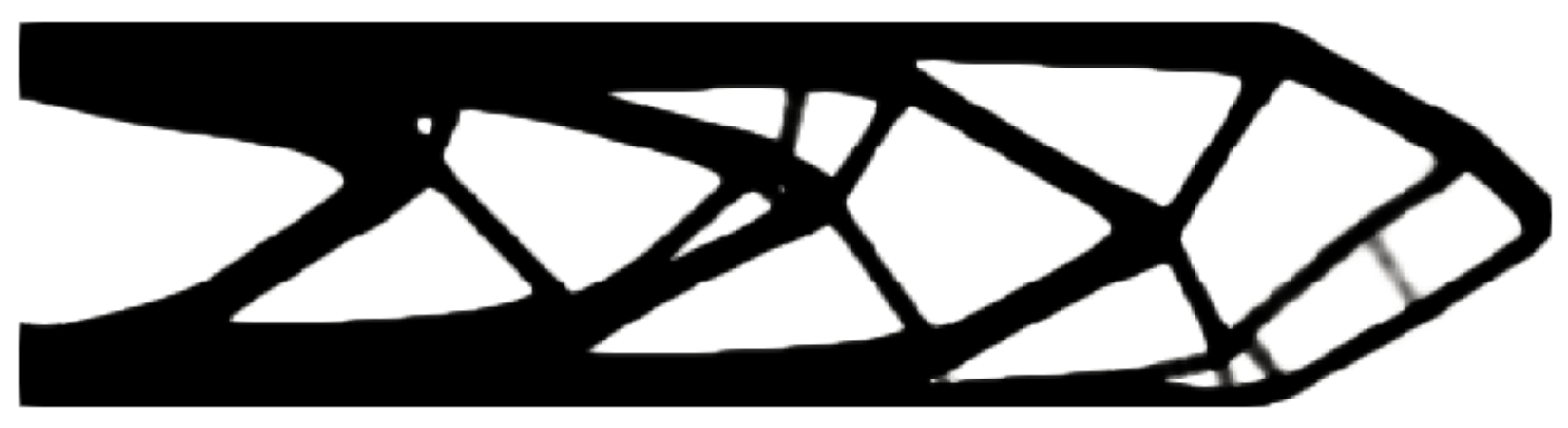} } & 
    \subfloat[$\bar{f} = 240$ kN, $\rho_{t} = 0.01$]{
    \includegraphics[width=0.32\linewidth]{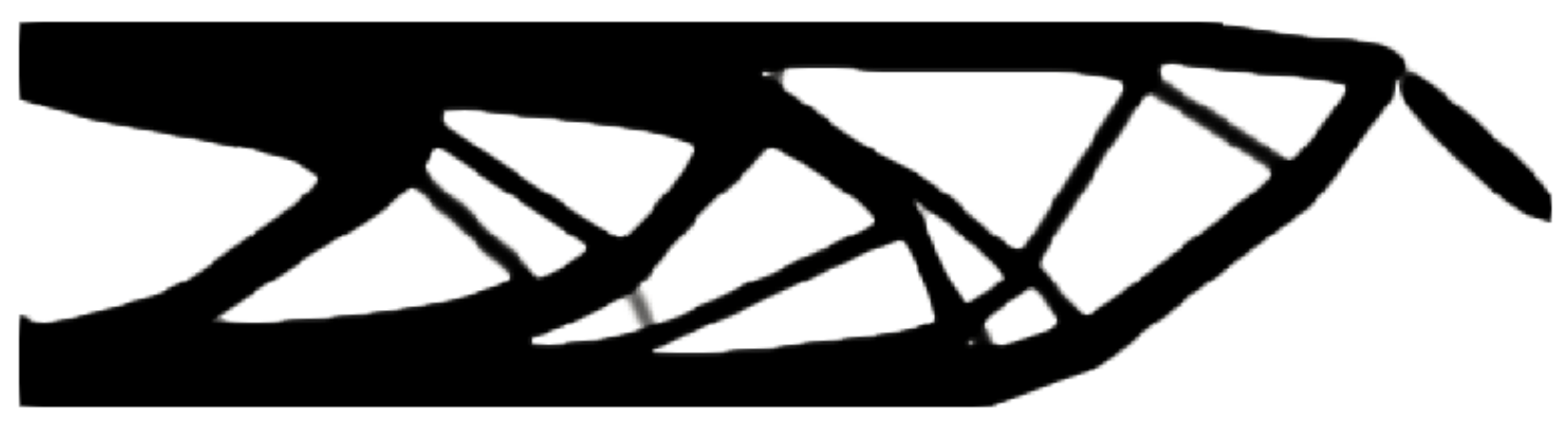} } \\
    \newline
    \subfloat[$\bar{f} = 12$ kN, $\rho_{t} = 0.1$]{
    \includegraphics[width=0.32\linewidth]{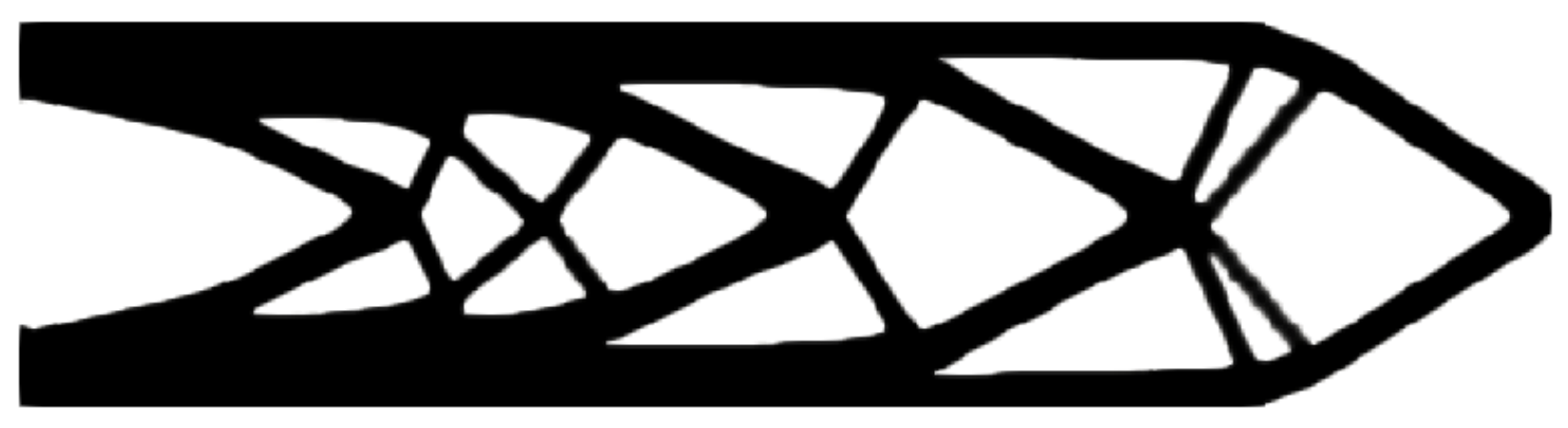} } &
    \subfloat[$\bar{f} = 144$ kN, $\rho_{t} = 0.1$]{
    \includegraphics[width=0.32\linewidth]{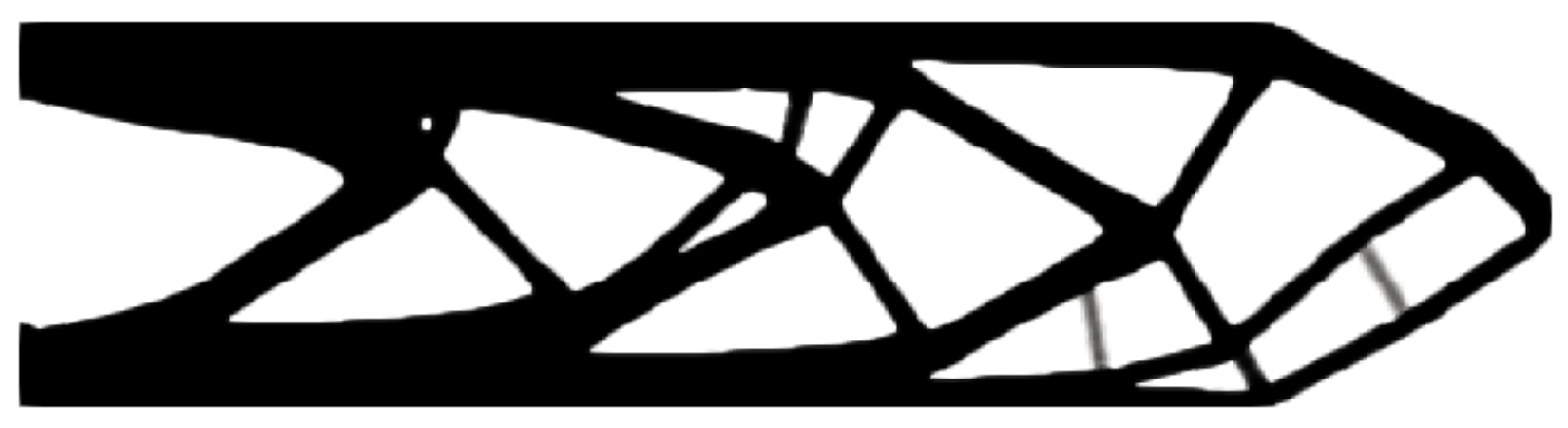} } & 
    \subfloat[$\bar{f} = 240$ kN, $\rho_{t} = 0.1$]{
    \includegraphics[width=0.32\linewidth]{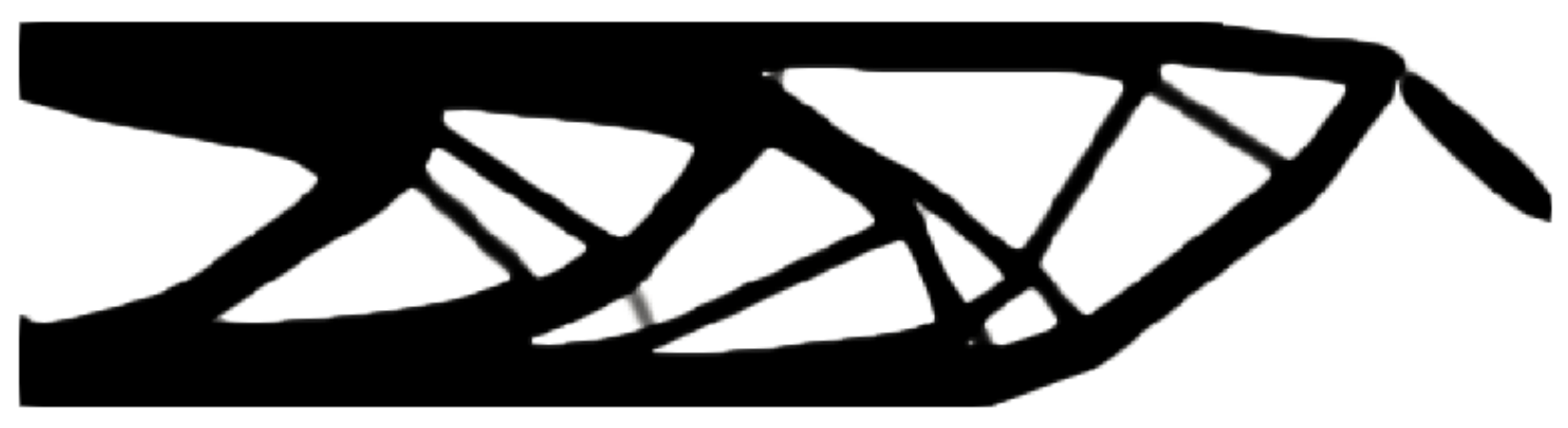} } \\
    \end{tabular}
\caption{Geometric nonlinear example: optimized designs obtained for increasing values of the applied point load (columns, left to right), and of the threshold density (rows, top to bottom) $\bar{f}$. For the highest value of the applied load ($\bar{f} = 240$ kN), the use of a threshold value $\rho_{t} < 0.01$ caused serious convergence issue in the nonlinear state solution, and therefore the designs are not reported.}
\label{fig:GeomNlnr_OptimizedDesignsAll}
\end{figure*}

The optimization problem is solved considering different values of $\bar{f}$ and $\rho_{t}$. The resulting optimized designs are shown in \autoref{fig:GeomNlnr_OptimizedDesignsAll}. For the design with applied load $\bar{f} = 12 \ \text{kN}$, the structural responses remain in the linear regime (for different threshold densities) and the optimized designs resemble the linear compliance responses given in \autoref{fig:linCompl_topResults}. However, as the applied load increases, the geometric nonlinearity becomes dominant, and the optimizer creates some pronounced features to minimize the compliance as can be seen in the designs with $\bar{f} = 240 \ \text{kN}$. We note that the results corresponding to $\rho_{t} = n.a.$ and $\bar{f} = 240 \ \text{kN}$ are excluded from \autoref{fig:GeomNlnr_OptimizedDesignsAll} because elements with a relative density very close to zero undergo excessive deformation resulting in convergence difficulties in the forward analysis. This instabilities are avoided in the present approach as low-density elements are removed from the forward analysis. This strategy is different from the approach given in \cite{BPS:00}, where the corresponding \ac{DOFs} for low relative density elements are removed from the forward analysis convergence criteria.

\begin{figure*}[!ht]
 \centering
  \subfloat[(a)]{
   \includegraphics[width = 0.30\linewidth]
    {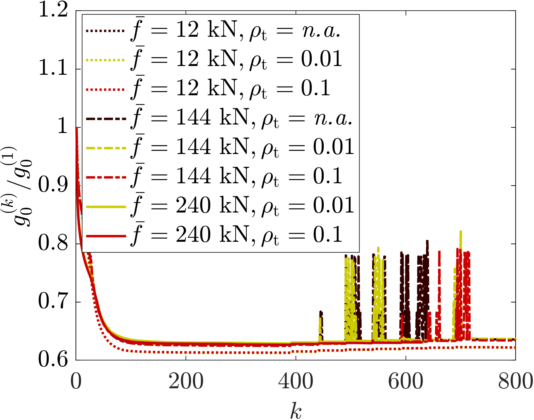}
   }
   \quad
  \subfloat[(b)]{
   \includegraphics[width = 0.30\linewidth]
    {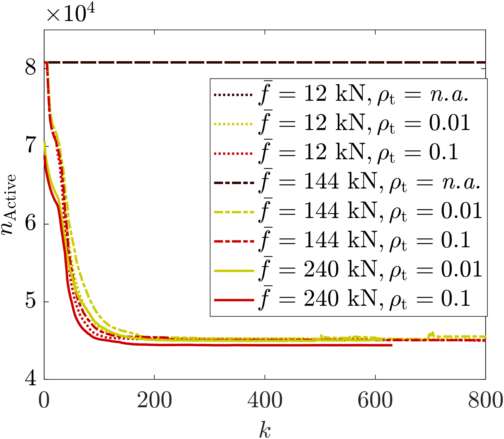}
   }
   \quad
  \subfloat[(c)]{
   \includegraphics[width = 0.30\linewidth]
    {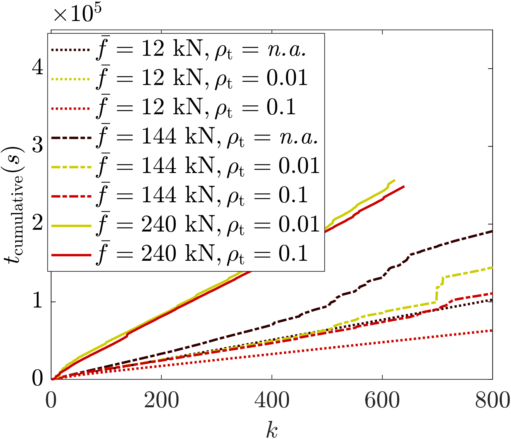}
   }
  \caption{Optimization history for the geometric nonlinear problem: (a) the evolution of the normalized compliance ($g^{(k)}_{0}/g^{(1)}_{0}$, (b) the number of \enquote{active} \ac{DOFs} in the state analysis ($n_{\rm Active}$), and (c) the cumulative computational time.}
 \label{fig:geomNlnr_2D_optResults}
\end{figure*}

The evolution of the objective function is shown in \autoref{fig:geomNlnr_2D_optResults}(a). The objective values are normalized with respect to their initial value. It can be seen that for some applied load and threshold density cases, there are oscillations in the compliance of the system that are caused due to sudden appearing/disappearing of some features in the topology optimized designs as $\eta$ increases. This is related to the system's highly nonlinear behavior and redistribution of the strain energy as the \ac{SIMP} exponent, $\eta$, increases.

The evolution of the number of \enquote{active} \ac{DOFs} in the state analysis and the cumulative computational time (i.e., \autoref{fig:geomNlnr_2D_optResults}(b, c)) clearly show the performance of the developed elemental removal strategy in reducing the computational cost in nonlinear structural systems by about 50\%. \autoref{fig:geomNlnr_2D_optResults}(c) clearly illustrates an increase in the cumulative computational time as the nonlinearity increases. This is mainly due to an increase in the number of forward analysis iterations (i.e., Newton-Raphson iterations) as the applied load increases, and the system becomes highly nonlinear.

\subsection{Topology Optimization with eigenvalues}
 \label{sSec:vibrationsAndbuckling}

Eigenvalue criteria play a crucial role for many design problems, as several important response quantities can be characterized as the minimizer of a Rayleigh quotient $\mathcal{R}(\boldsymbol{u}(\boldsymbol{\phi}))$ \cite{Washizu:75}, i.e.,
\begin{equation}
 \label{eq:RayMinimumCharacterization}
  \lambda_{1}(\boldsymbol{\phi}) = \min_{i=1,\ldots,n}
  \left\{
  \min_{\boldsymbol{v}_{i}\in\mathbb{R}-\boldsymbol{0}}
  \mathcal{R}( \boldsymbol{v}_{i} ) :=
  \frac{\boldsymbol{v}^{T}_{i}A(\boldsymbol{\phi})\boldsymbol{v}_{i}}
  {\boldsymbol{v}^{T}_{i}B(\boldsymbol{\phi})\boldsymbol{v}_{i}}
  \right\},
\end{equation}
where $A(\boldsymbol{\phi})$ and $B(\boldsymbol{\phi})$ discretize some energy-related operators (strain, kinetic, stress energies). For such cases, the state problem, derived from the condition $\partial_{\boldsymbol{v}} \mathcal{R}(\boldsymbol{u}(\boldsymbol{\phi})) = \mathbf{0}$ takes the form of a \ac{GEP}.

The following two sections present two popular examples of topology optimization for eigenvalues, i.e., the maximization of the fundamental frequency of vibration and the fundamental buckling load factor. Specifically, with these examples we challenge the element removal strategy to check (1) the effective reduction of the computational cost, (2) the robustness with respect to the appearance of artificial modes, and (3) the ability to reintroduce material when needed.

We also want to test the sensitivity of the element removal strategy with respect to the artificial modes phenomenon, giving some guidance about the selection of the threshold density to avoid its triggering. In order to identify artificial (buckling or vibration) modes we adopt the measure proposed by \cite{GM:15}, based on the ratio between the strain energy on regions with low relative density, and the overall one. Calling $\boldsymbol{v}_{i, \rm low}$ the subset of displacements on the elements with $\rho_{e} < 0.1$, the $i$-th mode will be regarded as artificial if
\begin{equation}
 \label{eq:strainEnergyRatioSPuriousModes}
  \Phi_{\rm rat} = \frac{\boldsymbol{v}^{T}_{i, \rm low}
  K(\boldsymbol{\phi})\boldsymbol{v}_{i}}
  {\boldsymbol{v}^{T}_{i}K(\boldsymbol{\phi})\boldsymbol{v}_{i}}
  > 0.5.
\end{equation}

The ability to reintroduce material is now particularly important, since the relationship $\boldsymbol{\phi}\mapsto \lambda_{1}(\boldsymbol{\phi})$ is generally non-monotonous, and we have to expect a more complicated design evolution. Moreover, reinforcement problems are often encountered in buckling design, and for these we precisely require the optimizer to reallocate the material starting from a given, well-defined, initial configuration.

\begin{figure*}[tb]
 \centering
  \subfloat[(a)]{
   \includegraphics[width=0.45\linewidth]
    {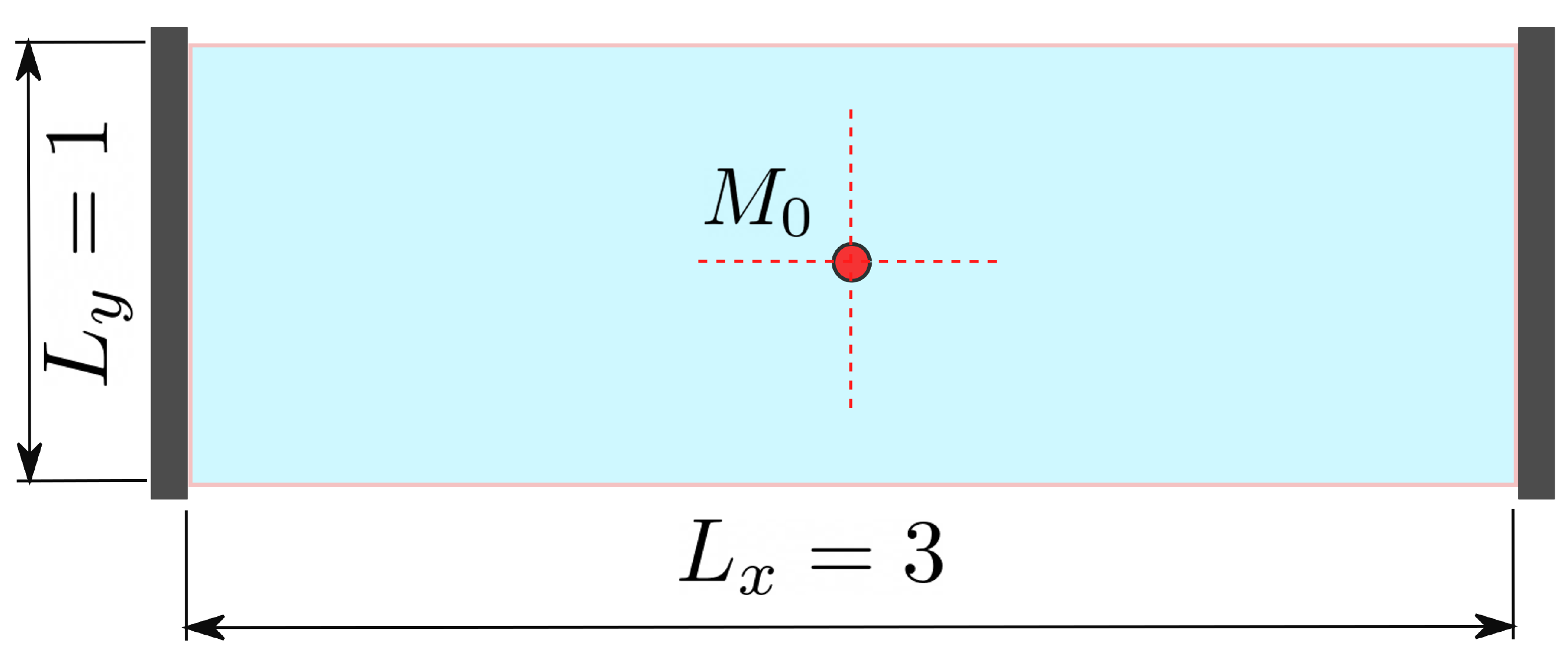}
   }
  \quad
  \subfloat[(b) $\rho_{t} = n.a.$, $\omega_{1} = 0.3045$]{
   \includegraphics[width=0.45\linewidth]
    {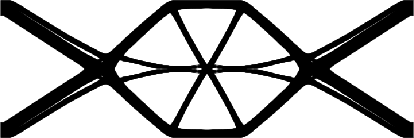}
   }
  \\
  \subfloat[(c) $\rho_{t} = 0.05$, $\omega_{1} = 0.3022$]{
   \includegraphics[width=0.45\linewidth]
    {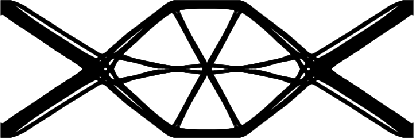}
   }
  \quad
  \subfloat[(d) $\rho_{t} = 0.1^{\ast}$, $\omega_{1} = 0.3038$]{
   \includegraphics[width=0.45\linewidth]
    {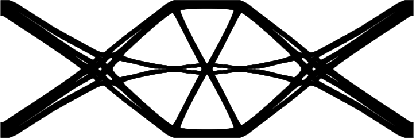}
   }
  \caption{Frequency maximization example: geometrical setup (a), reference design (b) and designs obtained for increasing values of the removal threshold (c,d). The threshold value $\rho_{t} = 0.1^{\ast}$ is obtained starting with $\rho_{t} = 0$ and increasing the threshold by $\Delta \rho_{t} = 0.02$ each 50 re-design steps.}
 \label{fig:vib_clamped_RAMP_topComparison}
\end{figure*}

\subsubsection{Maximization of the fundamental frequency of vibration}
 \label{ssSec:maxVibration}

We refer to the clamped beam sketched in \autoref{fig:vib_clamped_RAMP_topComparison}(a), discretized by $720\times 240$ isoparametric $\mathcal{Q}_{4}$ elements, for nearly $3.5  \times 10^{5}$ \ac{DOFs}. The magnitude of the nonstructural mass $|M_{0}|$, placed at the center of the design domain, and active on both \ac{DOFs}, is $15\%$ of the maximum allowed structural mass corresponding to $V_{\rm max} = 0.3$.

We aim at maximizing the fundamental frequency of vibration ($\omega_{1}$) \cite{MKC:95, Pedersen:00} and the optimization problem reads
\begin{equation}
 \label{eq:OptimizationProblem_maxFrequency}
  \begin{aligned}
   \min_{\boldsymbol{\phi}}
   & \: \Lambda^{-1}[ \omega_{i} ] \\
   \text{s.t.}
   & \: (K(\boldsymbol{\phi}) - \omega^{2}_{i} M(\boldsymbol{\phi}))\boldsymbol{v}_{i} = \boldsymbol{0}
   \: , \qquad \boldsymbol{v}_{i} \neq \boldsymbol{0} \\
   & V(\boldsymbol{\rho}) - V_{\rm max} \leq 0 \\
   & 0 \leq \phi_{i} \leq 1 \qquad i = 1, \ldots, m   
  \end{aligned},
\end{equation}
where $M(\boldsymbol{\phi})$ is the overall mass matrix (accounting for the nonstructural contribution $M_{0}$) and $\boldsymbol{v}_{i}$ are the vibration modes, normalized such that $\boldsymbol{v}^{T}_{i} M(\boldsymbol{\phi})\boldsymbol{v}_{j} = \delta_{ij}$.

At each optimization step we compute the lowest eight eigenpairs $(\omega^{2}_{i}, \boldsymbol{v}_{i})$ by solving the generalized eigenvalue problem associated with $K(\boldsymbol{\phi})$ and $M(\boldsymbol{\phi})$. Then, the inverses of the lowest $q = 4$ frequencies are aggregated in the Kreisselmeier-Steinhauser (KS) function \cite{KS:80}:
\begin{equation}
 \label{eq:KSaggregation-frequencies}
  \Lambda[\omega_{i}] = \omega^{-1}_{1} + \frac{1}{\alpha} \sum^{q}_{j=1} e^{\alpha (\omega^{-1}_{i} - \omega^{-1}_{1})},
\end{equation}
giving a lower bound to $\omega_{1}$. The sensitivity of \autoref{eq:KSaggregation-frequencies} is
\begin{equation}
 \label{eq:sensitivityFrequencyMaxExample}
  \frac{\partial \Lambda[\omega_{i}]}{\partial \rho_{e}} =
  \frac{1}
  {\sum^{q}_{i=1} e^{\alpha \omega^{-1}_{i}}}
  \sum^{q}_{i = 1} \frac{e^{\alpha \omega^{-1}_{i}}}
  {2\sqrt{\omega_{i}}}
  \boldsymbol{v}^{T}_{i}
  \left(
  \frac{\partial K}{\partial \rho_{e}} - \omega^{2}_{i}
  \frac{\partial M}{\partial \rho_{e}}
  \right)
  \boldsymbol{v}_{i}.
\end{equation}

The structural mass is parametrized by the linear interpolation $\varrho(\rho_{e}) = \rho_{e}\varrho_{0}$, where $\varrho_{0} = 1$, and the stiffness by the \ac{RAMP} interpolation \cite{SS:01}
\begin{equation}
 \label{eq:vibrations_RAMP}
  E(\rho_{e}) = E_{\rm min} + (E_{0}-E_{\rm min})
  \frac{\rho_{e}}
  {1 + \eta ( 1 - \rho_{e} )},
\end{equation}
starting with $\eta = 0$ and increasing the penalization by $\Delta\eta = 1$ at every 20 re-design steps, up to $\eta_{\rm max} = 24$. We set $E_{0} = 1$ and $E_{\rm min} = 10^{-6}$ for the case when $\rho_{t} = n.a.$, otherwise $E_{\rm min} = 0$. The filter size is $r_{\rm min} = 1.67 \times 10^{-2}L$ (4 elements) and the curvature regularization is set to $\beta = 40$. The aggregation parameter is set to $\alpha = 16$, and the optimization is run for 500 steps.

\begin{figure*}[tb]
 \centering
  \subfloat[(a) Reference design ($\rho_{t} = n.a.$)]{
   \includegraphics[width = 0.3\linewidth]
    {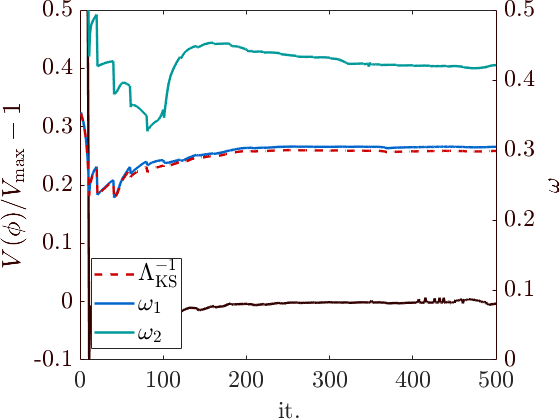}
   }
  \quad
  \subfloat[(b) $\rho_{t} = 0.02$]{
   \includegraphics[width = 0.3\linewidth]
    {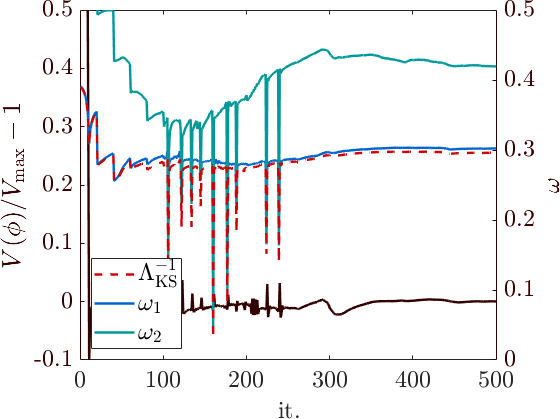}
   }
  \quad
  \subfloat[(c) $\rho_{t} = 0.05$]{
   \includegraphics[width = 0.3\linewidth]
    {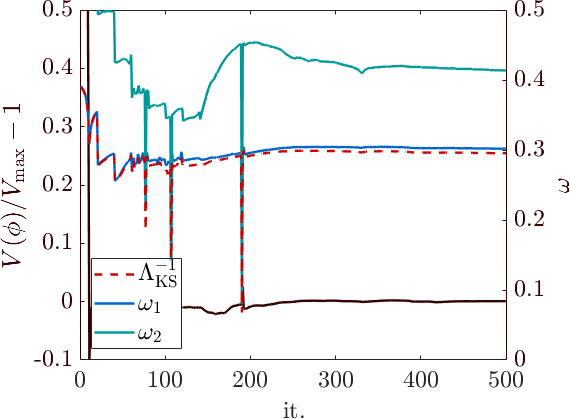}
   }
  \\
  \subfloat[(d) $\rho_{t} = 0.1^{\ast}$]{
   \includegraphics[width = 0.3\linewidth]
    {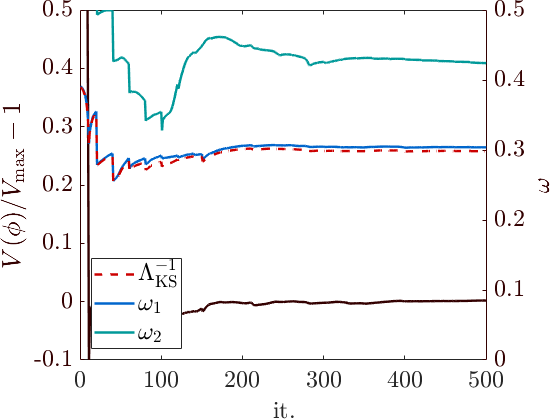}
   }
  \quad
  \subfloat[(e)]{
   \includegraphics[width = 0.28\linewidth]
    {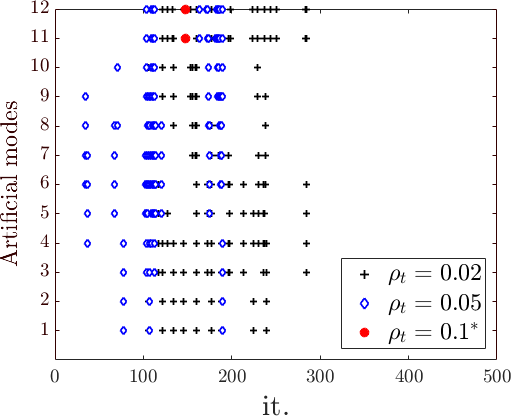}
   }
  \quad
  \subfloat[(f)]{
   \includegraphics[width = 0.3\linewidth]
    {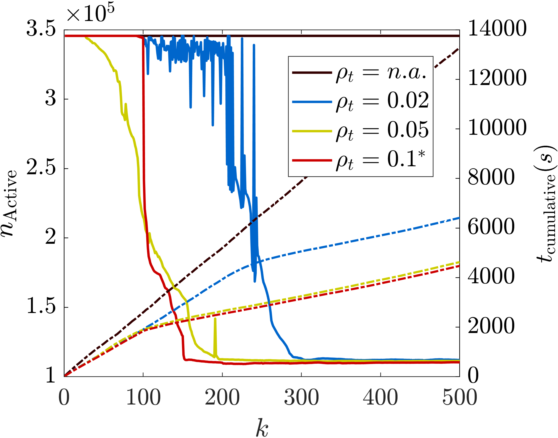}
   }
  \caption{Results for the frequency maximization example. The optimization histories (a-d) show the evolution of the volume constraint (black curve plotted against the left axis) and of the two lowest frequencies together with their KS aggregation function (plotted against the right axis), for several values of the removal threshold. We recall that the value $\rho_{t} = 0.1^{\ast}$ is reached by starting with $\rho_{t} = 0$ and increasing the threshold by $\Delta \rho_{t} = 0.02$ every 50 redesign steps. The marks in subplot (e) correspond to the occurrence of artificial vibration modes at a given iteration, based on the criterion of \autoref{eq:strainEnergyRatioSPuriousModes}. (f) shows the number of active \ac{DOFs} (continuous curves plotted against the left axis) and the cumulative CPU time (dashed curves plotted against the right axis).}
 \label{fig:vib_clamped_RAMP_optimizationResults}
\end{figure*}

The reference design is displayed in \autoref{fig:vib_clamped_RAMP_topComparison}(b) and the evolution of the lowest two frequencies, of their KS aggregation and of the volume constraint is shown in \autoref{fig:vib_clamped_RAMP_optimizationResults}(a). Except from a few jumps that correspond to $\eta$ increases in the early optimization steps, the optimization history is smooth and no artificial vibration modes are detected. The fundamental frequency attains the reference value $\omega_{1} = 0.3045$.

A comment on the selection of the threshold $\rho_{t}$ is now in order. On one hand, the dynamic nature of the objective (i.e., the competition between strain and kinetic energy) makes the grayscale regions more persistent \cite{ZWR+:16}. Thus, a very small magnitude of the threshold would not effectively remove the low relative density elements, slowing down the process and causing numerical issues. From our numerical studies, we observed reasonable results when choosing $\rho_{t} > 0.01$. On the other hand, the choice of a very large threshold value from the beginning of the optimization may hamper the development of the structural layout, promoting material disconnection.

We consider the threshold values $\rho_{t} = 0.02$, $0.05$ and $0.1$, and for the latter case we set $\rho_{t} = 0$ at the beginning of the optimization, then we increase it by $\Delta \rho_{t} = 0.02$ at every 50 optimization steps, up to $\rho_{t} = 0.1$. To highlight this, in the following we will denote this threshold value by $\rho_{t} = 0.1^{\ast}$. The corresponding optimization histories are shown in \autoref{fig:vib_clamped_RAMP_optimizationResults}(b-d), and the designs corresponding to $\rho_{t} = 0.05$ and $\rho_{t} = 0.1^{\ast}$ are in \autoref{fig:vib_clamped_RAMP_topComparison}(c,d).

Choosing $\rho_{t} = 0.02$ we obtain jumps in the frequencies at several steps in the first half of the optimization history. These are due to the appearance of artificial modes in the low relative density regions, where $\rho_{e}$ is slightly above $\rho_{t}$, and thus are not removed (we reiterate that we have set $E_{\rm min} = 0$ in \autoref{eq:vibrations_RAMP}). This unstable behavior can be cured by increasing the threshold density. For $\rho_{t} = 0.05$, artificial modes occur at three iterations only, whereas for the case $\rho_{t} = 0.1^{\ast}$  there are no artificial modes at all. It is remarkable to note that, despite some artificial modes and oscillations that may occur in the optimization history, the final value attained by $\omega_{1}$ is nearly the same as the reference one for all the three threshold density values (relative differences are in the order of $0.9\%$, $0.7\%$ and $0.2\%$, respectively).

\begin{figure*}[tb]
 \centering
  \subfloat[(a) $\rho_{t} = 0.02$]{
   \includegraphics[width = 0.3\linewidth]
    {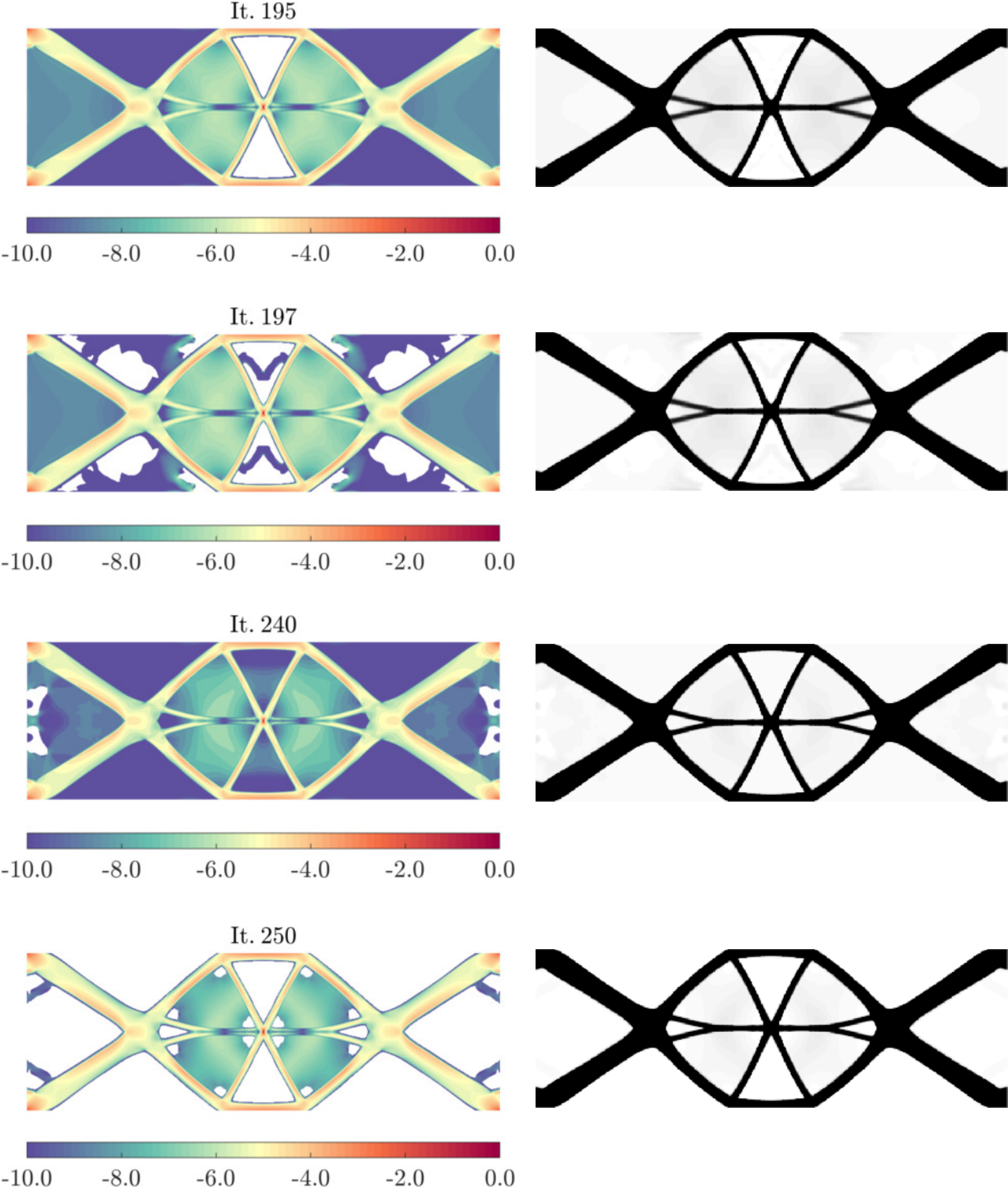}
  }
  \quad
  \subfloat[(b) $\rho_{t} = 0.05$]{
   \includegraphics[width = 0.3\linewidth]
    {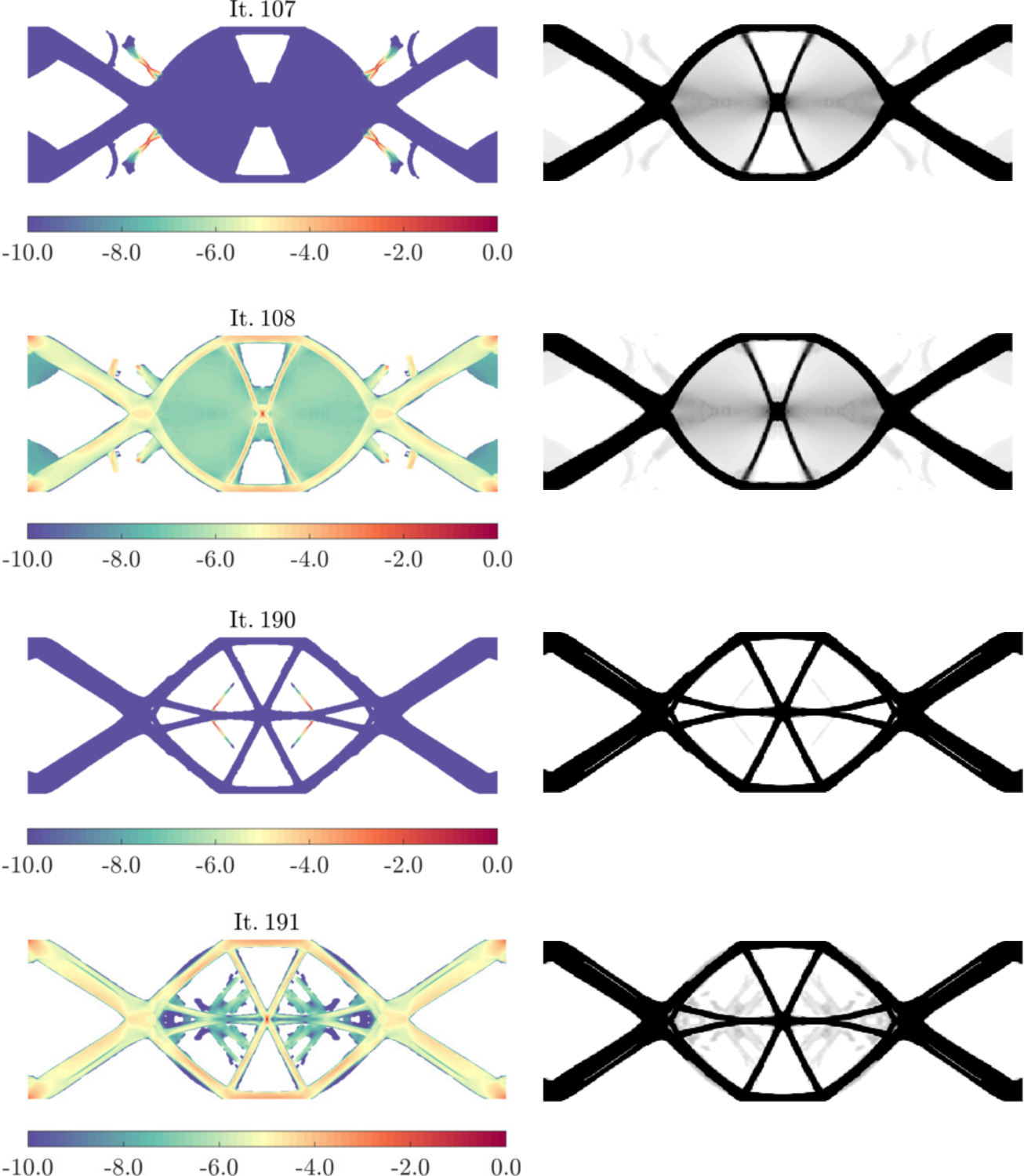}
  }
  \quad
  \subfloat[(c) $\rho_{t} = 0.1^{\ast}$]{
   \includegraphics[width = 0.3\linewidth]
    {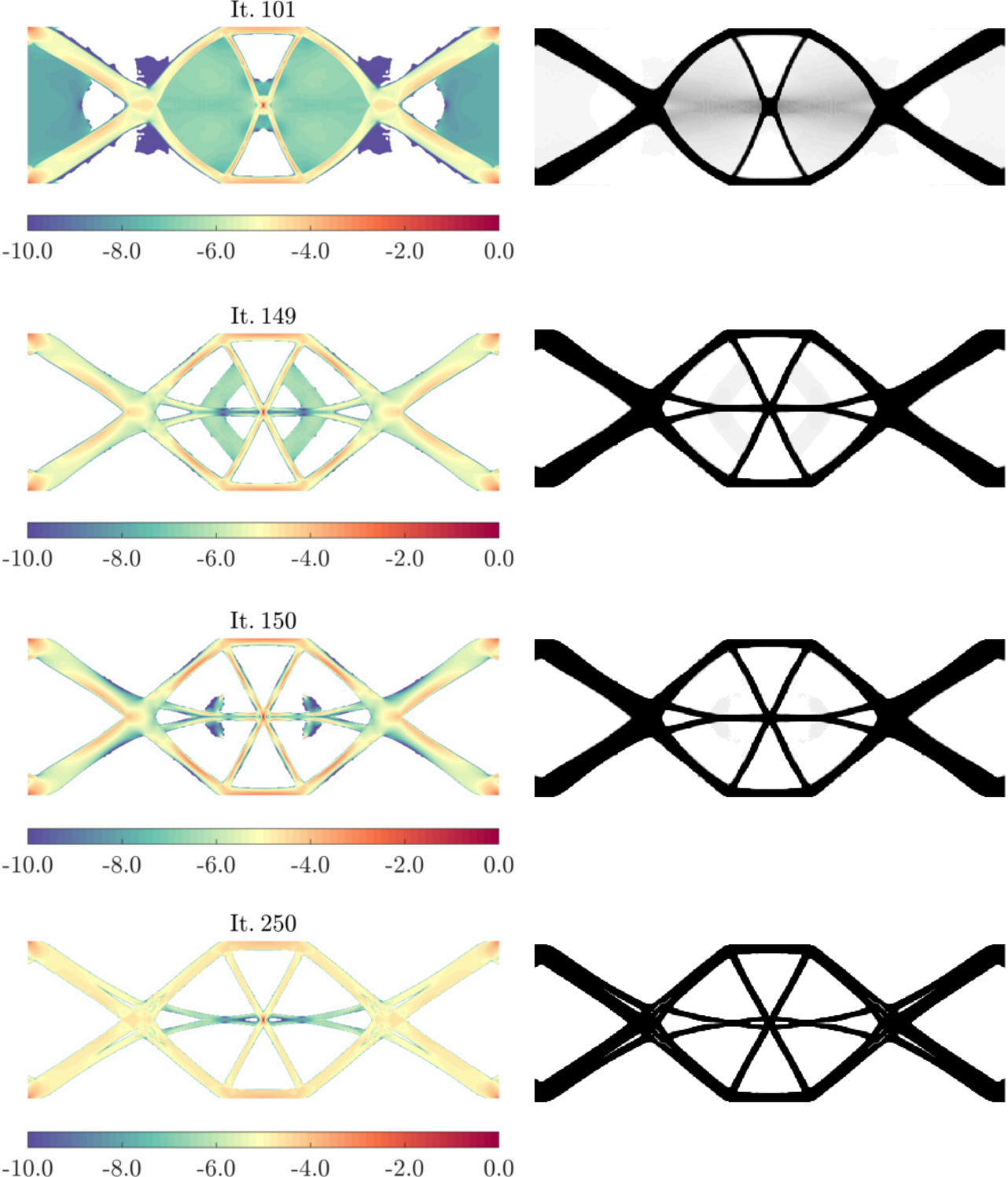}
  }
  \caption{Plots in the left columns show the logarithm of the strain energy density (normalized with respect to the maximum), corresponding to the fundamental vibration mode at some relevant optimization steps. Only the regions where $\rho_{e}\geq \rho_{t}$ are shown and the ``white'' regions are not modeled in the finite element analysis and replaced by fictitious boundary conditions. Plots in the right columns show the corresponding distribution of the volume fractions $\rho_{e}$.}
 \label{fig:vib_clamped_RAMP_SEDcontours}
\end{figure*}

To further discuss the evolution of the design and to grasp the origin of the discontinuities in the optimization history, we refer to \autoref{fig:vib_clamped_RAMP_SEDcontours}, showing the strain energy density distribution on the \enquote{active} elements at some relevant iterations.

For the case $\rho_{t} = 0.02$ (see \autoref{fig:vib_clamped_RAMP_SEDcontours}(a)), although the design shows large regions with almost white material, we still have $\rho_{e} > \rho_{t}$ in these regions. Therefore, the element removal is slow and the material reintroduction is also very frequent (see iterations 195-197 and \autoref{fig:vib_clamped_RAMP_optimizationResults}(f), showing the number of active elements). The persistence of these large regions with low stiffness to mass ratio clearly promotes the appearance of artificial modes. However, the structural configuration becomes well defined after iteration 250, and the optimization proceeds by slowly removing all the remaining gray regions.

\autoref{fig:vib_clamped_RAMP_SEDcontours}(b), that corresponds to the threshold $\rho_{t} = 0.05$, clearly illustrates the origin of a spurious mode and its immediate removal between iterations 107-108 and 190-191. In both cases, some low density material, that is slowly fading away but is still above the threshold density, is shaped in some bars surrounded by voids. Therefore, at iteration 107 (viz. 190), these bars vibrate with a very low frequency, originating a spurious mode. In the next iteration, these bars are made thicker and shorter (see iteration 108), or some surrounding material reappear to strengthen them (see iteration 191), causing a sudden increase in the frequency, as the vibration involves the whole structure again. We emphasize that one can distinctly observe the element reintroduction by looking at \autoref{fig:vib_clamped_RAMP_optimizationResults}(f).

For the case with $\rho_{t} = 0.1^{\ast}$, the evolution of the design and \enquote{active} elements is very smooth (see \autoref{fig:vib_clamped_RAMP_SEDcontours}(c)). Here, we point out that essentially all elements are modeled up to iteration 100, when the threshold density is increased to $\rho_{t} = 0.02$ and all the low relative density regions outside of the structural boundaries are removed at once (see also \autoref{fig:vib_clamped_RAMP_optimizationResults}(f) for the \enquote{active} element decrease). The same happens across iteration $150$, with the further increase of $\rho_{t}$, and after this point the structural configuration exhibit only minor changes.

Finally, we notice the huge computational savings associated with the element removal approach by looking at \autoref{fig:vib_clamped_RAMP_optimizationResults}. The overall CPU time is cut by $50\%$ for the low value $\rho_{t} = 0.02$ and by $70\%$ for the larger values $\rho_{t} = 0.05$ and $\rho_{t} = 0.1^{\ast}$.

\begin{figure*}[tb]
 \centering
  \subfloat[(a) Geometrical setup and min. volume design]{
   \includegraphics[width=0.45\linewidth]
    {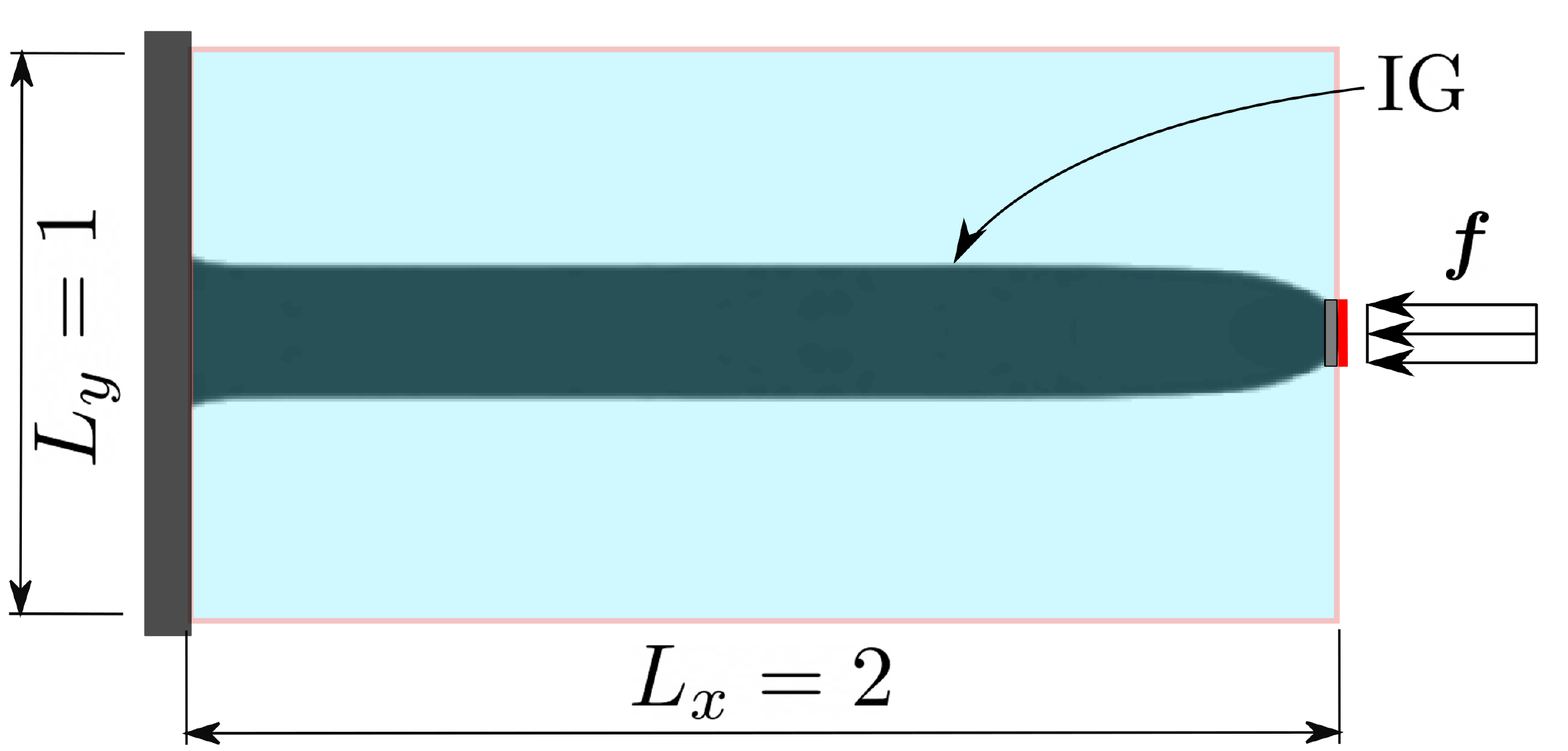}
   }
  \\
  \subfloat[(b) $\rho_{t} = n.a.$, $\lambda_{1} = 5.491$]{
   \includegraphics[width=0.40\linewidth]
    {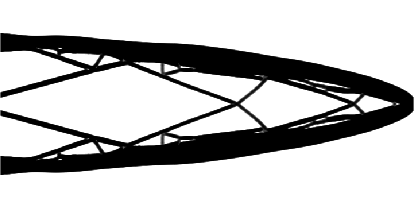}
   }
  \subfloat[(c) $\rho_{t} = 0.1^{\ast}$, $\lambda_{1} = 5.596$]{
   \includegraphics[width=0.40\linewidth]
    {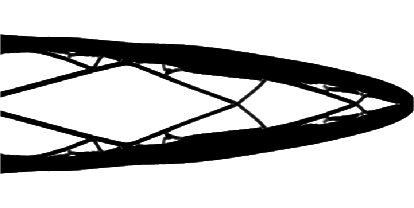}
   }
  \\
  \subfloat[(d) $\rho_{t} = n.a.$, \ac{IG}, $\lambda_{1} = 4.621$]{
   \includegraphics[width=0.40\linewidth]
    {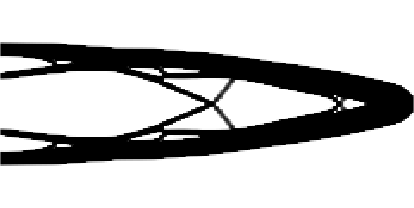}
   }
  \subfloat[(e) $\rho_{t} = 0.01$, \ac{IG}, $\lambda_{1} = 4.887$]{
   \includegraphics[width=0.40\linewidth]
    {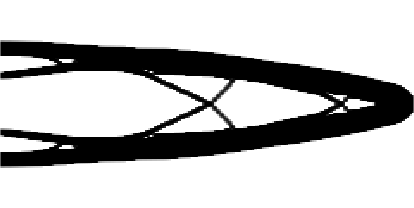}
   }
  \caption{Buckling load maximization example: (a) geometrical setup with superimposed the minimum volume design ($V(\boldsymbol{\phi}) \approx 0.22$) corresponding to the maximum allowed compliance $c_{\rm max} = 2.5 c_{\rm full}$. (b,c) are the designs obtained starting from the uniform initial guess $\boldsymbol{\phi} = 0.25$ and with $\rho_{t} = n.a.$ and $\rho_{t} = 0.1^{\ast}$, respectively. Designs in (d,e) are obtained starting from the minimum volume design, denoted as \ac{IG}, and with thresholds $\rho_{t} = n.a.$ and $\rho_{t} = 0.01$.}
 \label{fig:buckl_column_topComparison}
\end{figure*}

\subsubsection{Maximization of the Buckling Load Factor}
 \label{ssSec:maxBLF}

We consider the column of \autoref{fig:buckl_column_topComparison}(a) and the problem of maximizing its fundamental \ac{BLF} associated with the external load $\boldsymbol{f}$ \cite{GM:15,FS:19}, given the maximum volume fraction $V_{\rm max} = 0.25$. We also prescribe the upper bound on the pre-buckling compliance $c_{\rm max} = 2.5 c_{\rm full}$, being $c_{\rm full}$ the compliance of the fully solid column.

The domain is discretized by $720\times 360$ $\mathcal{Q}_{4}$ elements, for about $5.2  \times 10^{5}$ \ac{DOFs} and the compressive force is applied over the length $b=L_{y}/90$, with total magnitude $|\boldsymbol{f}| = 10^{-3}$.

The optimization problem reads
\begin{equation}
 \label{eq:OptimizationProblem_maxBLF}
  \begin{aligned}
   \min_{\boldsymbol{\phi}}
   & \: \Lambda^{-1}[ \lambda_{i} ] \\
   \text{s.t.}
   & \: (K(\boldsymbol{\phi})+\lambda G(\boldsymbol{\phi}, \boldsymbol{u}(\boldsymbol{\phi})))\boldsymbol{v}_{i} = \boldsymbol{0}
   \: , \qquad \boldsymbol{v}_{i} \neq \boldsymbol{0} \\
   & K(\boldsymbol{\phi})\boldsymbol{u}(\boldsymbol{\phi}) = \boldsymbol{f} \\
   & \boldsymbol{u}(\boldsymbol{\phi})^{T}\boldsymbol{f} - c_{\rm max} \leq 0 \\
   & V(\boldsymbol{\rho}) - V_{\rm max} \leq 0 \\
   & 0 \leq \phi_{i} \leq 1 \qquad i = 1, \ldots, m   
  \end{aligned},
\end{equation}
where $G(\boldsymbol{\phi}, \boldsymbol{u}(\boldsymbol{\phi}))$ is the stress stiffness matrix defined according to the linearized buckling theory \cite{CMP+:01} and $\boldsymbol{v}_{i}$ represents the buckling modes, normalized such that $\boldsymbol{v}^{T}_{i} G(\boldsymbol{\phi},\boldsymbol{u}(\boldsymbol{\phi}))\boldsymbol{v}_{j} = -\delta_{ij}$.

At each re-design step, we compute the lowest 16 eigenpairs $(\lambda_{i}, \boldsymbol{v}_{i})$ by solving the generalized eigenvalue problem and the objective function $\Lambda[\lambda_{i}]$ is formed as in \autoref{eq:KSaggregation-frequencies}, by considering the inverses of the lowest $q = 12$ \acp{BLF}. Thus, the objective provides a lower bound to the \ac{BLF}, and its sensitivity reads
\begin{equation}
 \label{eq:sensitivityBucklingMaxExample}
  \frac{\partial \Lambda[\lambda_{i}]}{\partial \rho_{e}} =
  \frac{1}
  {\sum^{q}_{i=1} e^{\alpha \lambda^{-1}_{i}}}
  \sum^{q}_{i = 1} \frac{e^{\alpha \lambda^{-1}_{i}}}
  {2\sqrt{\lambda_{i}}}
  \boldsymbol{v}^{T}_{i}
  \left(
  \frac{\partial K}{\partial \rho_{e}} + \lambda_{i}
  \frac{\partial G}{\partial \rho_{e}} - \boldsymbol{a}^{T}_{i}
  \frac{\partial K}{\partial \rho_{e}}\boldsymbol{u}
  \right)
  \boldsymbol{v}_{i},
\end{equation}
where $\boldsymbol{a}_{i}$ are the adjoint vectors, one per buckling mode, solving the problem $K(\boldsymbol{\phi})\boldsymbol{a}_{i} = \boldsymbol{v}^{T}_{i}\partial_{\boldsymbol{u}}G(\boldsymbol{\phi},\boldsymbol{u}(\boldsymbol{\phi}))\boldsymbol{v}_{i}$ \cite{RGB:95}.

The interpolation functions adopted for the stiffness and stresses are \cite{BS:04}
\begin{equation}
 \label{eq:buckling_SIMP2SIMP}
  E(\rho_{e}) = E_{\rm min} + (E_{0}-E_{\rm min}) \rho^{\eta}_{e}
  \: , \qquad
  \sigma(\rho_{e}) = E_{0}\rho^{\eta}_{e},
\end{equation}
where $E_{0} = 1$, $E_{\rm min} = 10^{-6}E_{0}$ for the case when $\rho_{t} = n.a.$, otherwise $E_{\rm min} = 0$. We set $\eta = 3$ at the beginning of the optimization, then, starting from iteration 100, we increase the penalization by $\Delta \eta = 0.2$ at every 15 iterations, up to the value $\eta_{\rm max} = 6$. The filter radius is $r_{\rm min} = 1.1  \times 10^{-2}L_{y}$ (4 elements) and the curvature regularization parameter begins with $\beta = 4$ and is increased by $\Delta\beta = 2$ at every 50 iterations, up to $\beta_{\rm max} =24$, starting from iteration 350. We set the aggregation parameter to $\alpha = 16$ and the optimization is run for 750 iterations.

\begin{figure*}[tb]
 \centering
  \subfloat[(a) Reference design ($\rho_{t} = n.a.$)]{
   \includegraphics[width = 0.3\linewidth]
    {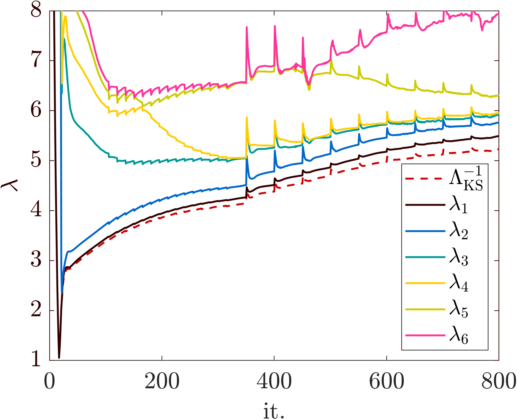}
   }
  \quad
  \subfloat[(b) $\rho_{t} = 0.01$]{
   \includegraphics[width = 0.3\linewidth]
    {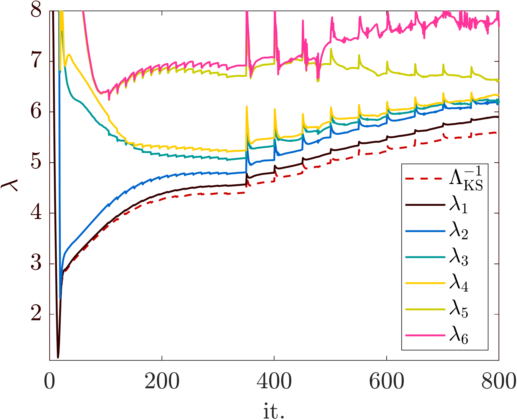}
   }
  \quad
  \subfloat[(c) $\rho_{t} = 0.1^{\ast}$]{
   \includegraphics[width = 0.3\linewidth]
    {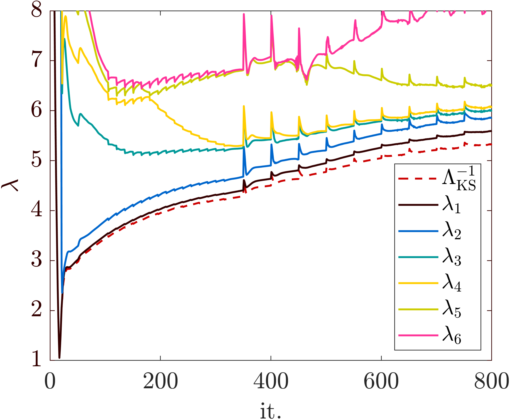}
   }
  \\
  \subfloat[(d) Start from \ac{IG} $\rho_{t} = n.a.$]{
   \includegraphics[width = 0.3\linewidth]
    {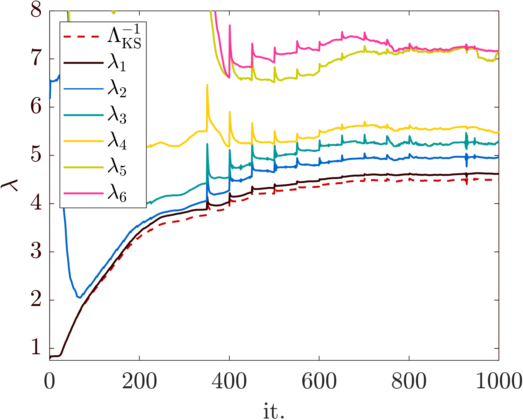}
   }
  \quad
  \subfloat[(e) Start from \ac{IG} $\rho_{t} = 0.1^{\ast}$]{
   \includegraphics[width = 0.3\linewidth]
    {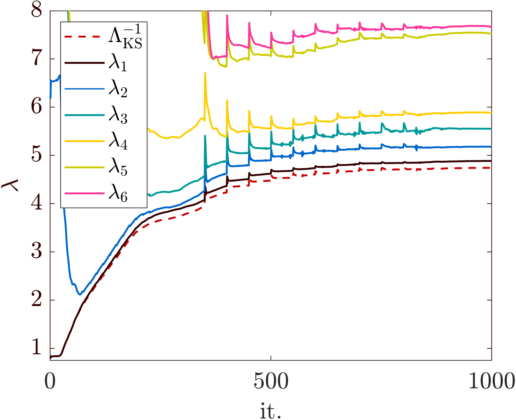}
   }
  \quad
  \subfloat[(f)]{
   \includegraphics[width = 0.3\linewidth]
    {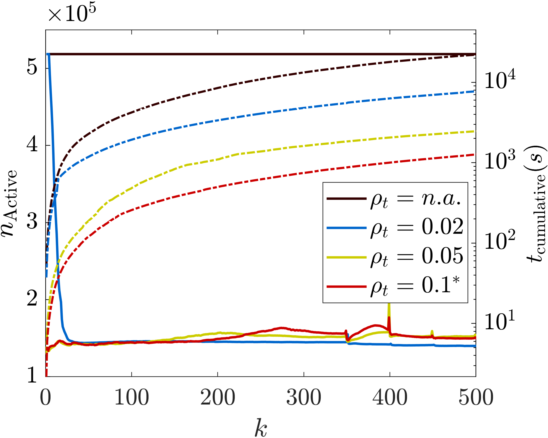}
   }
  \caption{Results for the buckling load maximization example. Plots (a-c) show the evolution of the lowest six buckling modes, together with the lower bound provided by the objective function $\Lambda^{-1}$, when starting from the uniform initial guess $\phi_{e}=0.25$. Plots (d,e) are referred to the buckling maximization problem starting from the minimum volume design (reinforcement problem). The number of active elements and progressive CPU time (in log scale) is shown in (f).}
 \label{fig:buckl_column_optimizationResults}
\end{figure*}

The reference design, obtained starting from the uniform material distribution $\phi_{e} = 0.25$, $e=1,\ldots, m$ is shown in \autoref{fig:buckl_column_topComparison}(b), whereas \autoref{fig:buckl_column_topComparison}(c) shows the design obtained with the element removal strategy for the threshold density value $\rho_{t} = 0.1^{\ast}$, reached by continuation as explained in \autoref{ssSec:maxVibration}. The evolution of the objective function and the lowest six \acp{BLF} for the standard approach, $\rho_{t} = 0.01$ and $\rho_{t} = 0.1$ are shown in \autoref{fig:buckl_column_optimizationResults}(a,b). The fundamental \ac{BLF} reaches the value $\lambda_{1} = 5.491$ for the reference design, $\lambda_{1} = 5.745$ for the design obtained with $\rho_{t} = 0.01$ and $\lambda_{1} = 5.901$ for the design obtained with $\rho_{t} = 0.1^{\ast}$. Again, the slightly higher value attained by the \ac{BLF} when using the element removal strategy is due to the suppression of low relative density elements around the boundaries. Similar considerations as made in \autoref{ssSec:maxVibration} for the frequency maximization hold, concerning the selection of the threshold density, and the use of a very low value excessively slows down the element removal. However, compared to the frequency optimization example, the process is now more robust with respect to low values of the removal threshold $\rho_{t}$ and, according to the criterion \autoref{eq:strainEnergyRatioSPuriousModes}, we have never encountered artificial buckling modes (see \autoref{fig:buckl_column_optimizationResults}(c)). Also, the grayscales seem less persistent than before. Both of these effects are due to the fact that the optimization is now indirectly driven also by the linear pre-buckling problem, since the stress stiffness operator $G(\boldsymbol{\phi},\boldsymbol{u}(\boldsymbol{\phi}))$ depends on the linear equilibrium displacement $\boldsymbol{u}$, and the stresses are further penalized by \autoref{eq:buckling_SIMP2SIMP}.

\begin{figure*}[!tb]
 \centering
 \subfloat[(a)]{
  \includegraphics[width = 0.11\linewidth]
   {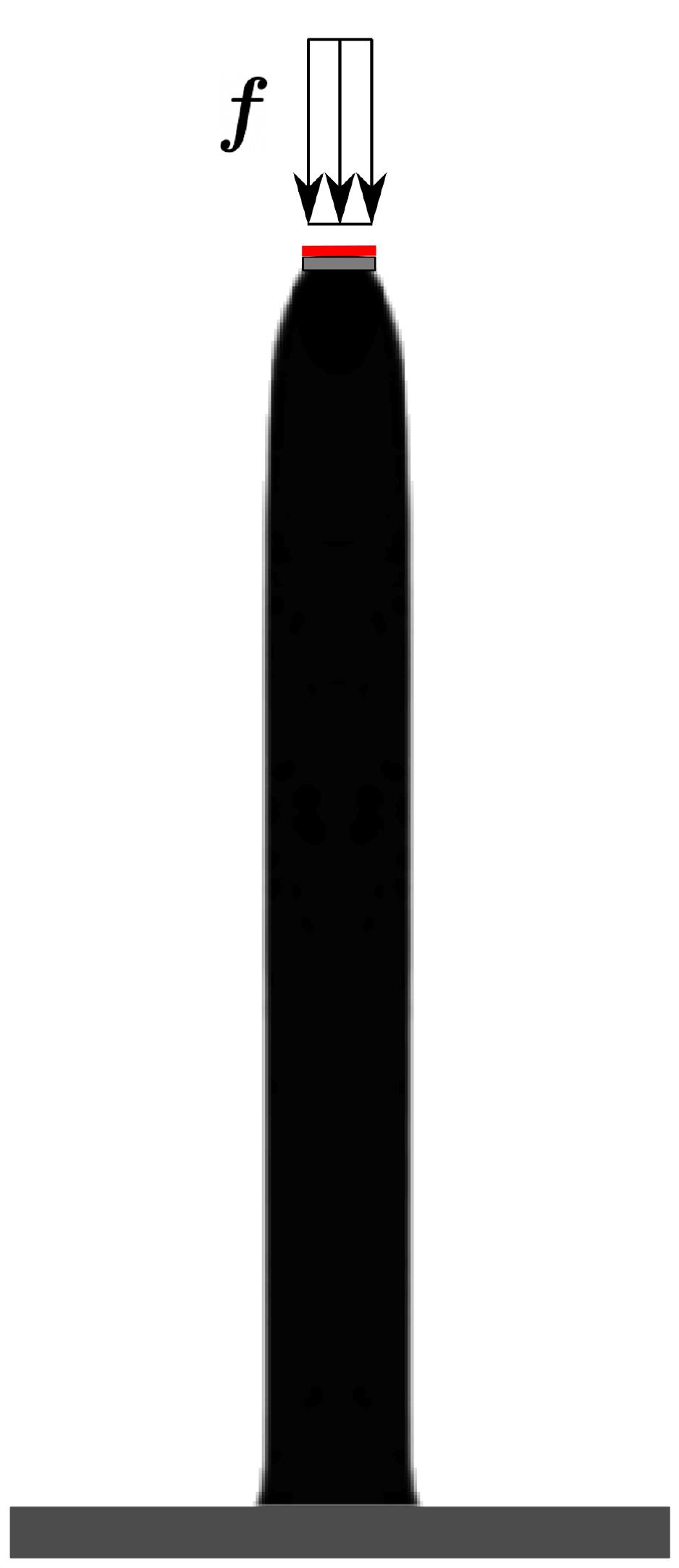}} \qquad
 \subfloat[(b)]{
  \includegraphics[width = 0.75\linewidth]
   {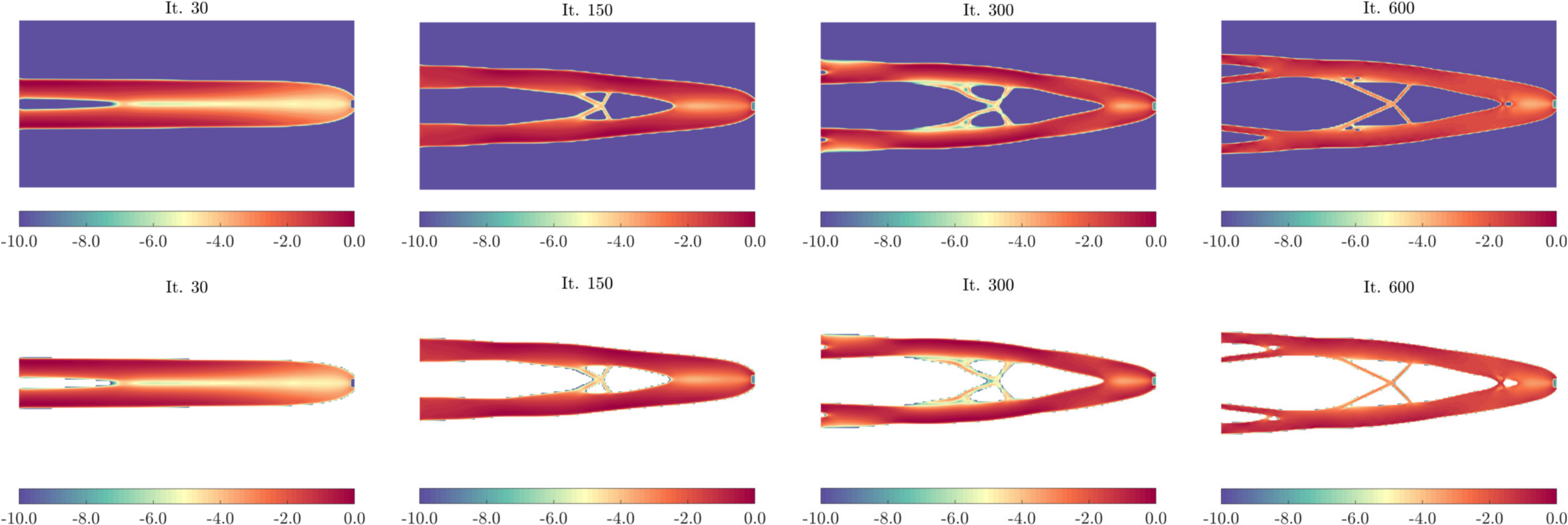}}
  \caption{Compressive column design problem: (a) minimum volume initial guess and (b) evolution of the design and of the strain energy density associated with the fundamental buckling mode, for the BLF maximization problem. Plots in the top row of (b) correspond to the standard approach ($\rho_{t} = n.a.$), while those in the bottom row correspond to the element removal approach with $\rho_{t} = 0.01$, and only the active elements are displayed.}
 \label{fig:buckl_column_SEDcontours}
\end{figure*}

Material reintroduction is not precluded by the element removal strategy. To highlight this, as it is one of the strong points of the present approach, we solve again problem \autoref{eq:OptimizationProblem_maxBLF} starting from a different initial guess. We consider the minimum volume design corresponding to the maximum allowed compliance $c_{\rm max} = 2.5 c_{\rm full}$, attaining the volume fraction $V(\boldsymbol{\phi})\approx 0.22$ (see \autoref{fig:buckl_column_SEDcontours}(a)). This thin bar is clearly prone to lateral buckling, and having $\lambda_{1} = 0.79$ it would buckle under the applied load. Increasing the \ac{BLF} clearly requires the material to spread away from the centerline, growing a nontrivial structural configuration, and this happens very effectively thanks to the material reintroduction facilitated by the \ac{HPM} (see \autoref{fig:buckl_column_SEDcontours}(b)). The reference design is shown in \autoref{fig:buckl_column_topComparison}(d) and the corresponding \ac{BLF} is $\lambda_{1} = 4.621$.  The design obtained with $\rho_{t} = 0.01$ is displayed in \autoref{fig:buckl_column_topComparison}(e) and, although it shows some minor difference in the configuration of the inner bars, we clearly see that it has essentially the same layout and an equivalent performance ($\lambda_{1} = 4.887$ is only $1\%$ lower than the reference value). \autoref{fig:buckl_column_optimizationResults}(d,e) show the objective and \acp{BLF} evolution for this case, and the similarity of the optimization histories is remarkable.

Finally, from \autoref{fig:buckl_column_optimizationResults}(f), showing the evolution in the number of active elements and the progressive CPU time, we appreciate the remarkable computational savings. The red and yellow continuous curves in subplot (f) show the number of active elements for the element removal approach and the corresponding number of elements where the relative density is above $\rho_{t} = 0.01$ for the standard approach, respectively, when starting from the minimum volume design of \autoref{fig:buckl_column_SEDcontours}(a). The evolution of there two curves is very similar, and further demonstrates the the ability of the proposed approach to reintroduce material, when needed.

%
\section{Conclusions
\label{sec:Conclusions}}
We have presented an adaptive element removal strategy, aimed at reducing the computational cost of density-based topology optimization. The method is based on the idea originally proposed in \cite{BT:03} of removing elements with a relative density below a given threshold from the state analysis, and to apply nodal boundary conditions to the \ac{DOF}s completely surrounded by removed elements. As a result, the size of the system of state equations solved at each re-design step and the corresponding computational effort is systematically reduced as the optimization progresses, as the solid layout emerges and many regions with low relative density appear.

A key feature of the proposed approach is that the removed elements can be readily reintroduced, as needed, since the entire design space remains active throughout the optimization process. The reintroduction of removed elements is driven purely by the formal sensitivity propagation effect due to the Heaviside projection, and the nonlinearity of this projection plays a key role in promoting the element reintroduction.

The effectiveness and convenience of the approach have been demonstrated on several examples from structural topology optimization, governed by both linear and nonlinear/eigenvalue state equations. We have also observed that the method helps to mitigate numerical instabilities associated with low relative density regions. First, the removal of elements with low relative density allows setting the stiffness lower bound to $E_{\rm min} = 0$ and prevents fictitious stiffening effects due to \enquote{void} for mechanisms design and nonlinear analyses. The detrimental effect of highly distorted elements in geometric nonlinear analysis, and artificial modes for eigenvalue-based problems can also be reduced by the proper selection of the element removal threshold.

We summarize the strengths of the presented method as follows:
\begin{itemize}
 \item It reduces the computational cost of the topology optimization procedure, without hampering the possibility to reintroduce material, if needed;
 \item It helps circumvent unrealistic effects and numerical instabilities associated with low relative density regions;
 \item It has promising generality and is easily extendable to topology optimization problems not considered herein as well as to more advanced preconditioned iterative solvers \cite{LHZ+:18}.
\end{itemize}

With regard to this last point, subsets of the current authors have also successfully implemented the presented approach for topology optimization problems governed by material nonlinearities, including plasticity \cite{Lotfi:14, CLG+:15}, and problems governed by fluid mechanics \cite{BRG:18}.

Finally, we again acknowledge that the action of introducing boundary conditions introduces a discontinuity in an otherwise continuous topology optimization problem. All of our work to date, including linear and nonlinear solid mechanics and fluid mechanics, has suggested this has little impact on the optimized design if $\rho_{t}$ is chosen reasonably ($\rho_{t} \leq 0.1$ herein). Nevertheless, it is certainly possible this could be problematic in examples yet considered by the authors. 

%
\section*{Acknowledgments
\label{sec:Acknowledgments}}
This material is based on research sponsored by the National Aeronautics and Space Administration (NASA) under Grant No. 80NSSC18K0428, Air Force Research Laboratory (AFRL) under agreement number FA8651-20-1-0007, and the National Science Foundation (NSF) under Grant No. 120103. The U.S. Government is authorized to reproduce and distribute reprints for Governmental purposes notwithstanding any copyright notation thereon. The views, findings, and conclusions or recommendations contained herein are those of the authors and should not be interpreted as necessarily representing the official policies or endorsements, either expressed or implied, of Johns Hopkins University or the U.S. Government.
%

%
\bibliography{behrou_to_nodalbc}

\end{document}